\documentstyle[12pt]{amsart}
\newcommand{\del}{\partial}
\newcommand{\bC}{{\Bbb C}}
\newcommand{\bP}{{\Bbb P}}
\newcommand{\bQ}{{\Bbb Q}}
\newcommand{\bT}{{\Bbb T}}
\newcommand{\bZ}{{\Bbb Z}}
\newcommand{\cF}{{\cal F}}
\newcommand{\cO}{{\cal O}}
\newcommand{\cP}{{\cal P}}
\newcommand{\cV}{{\cal V}}
\newcommand{\hS}{\widehat{S}}
\newcommand{\oL}{\overline{L}}
\newcommand{\op}{\overline{p}}
\newcommand{\orho}{\overline{\rho}}
\newcommand{\tC}{\widetilde{C}}
\newcommand{\tE}{\widetilde{E}}

\newcommand{\tq}{\widetilde{q}}
\newcommand{\tS}{\widetilde{S}}

\newcommand{\tZ}{\widetilde{Z}}
\newcommand{\lra}{\longrightarrow}
\newcommand{\inj}{\hookrightarrow}
\newcommand{\surj}{\hspace{3pt}\to \hspace{-19pt}{\rightarrow} \:\:}

\newtheorem{proposition}{Proposition}[section]
\newtheorem{lemma}[proposition]{Lemma}
\newtheorem{theorem}[proposition]{Theorem}
\newtheorem{definition}[proposition]{Definition}
\newtheorem{corollary}[proposition]{Corollary}
\newtheorem{remark}[proposition]{Remark}
\newtheorem{theorem0}{Theorem}
\newtheorem{lemma0}[theorem0]{Lemma}
\newtheorem{claim}[proposition]{Claim}

\begin{document}

\title[Cubic hypersurfaces]{A Prym construction for the cohomology of
a cubic hypersurface}
\author{E. Izadi}
\address{Department of Mathematics, Boyd Graduate Studies Research Center, 
University of Georgia, Athens, GA 30602-7403, USA}
\email{izadi@@math.uga.edu}

\begin{abstract}

Mumford defined a natural isomorphism between the intermediate
jacobian of a conic-bundle over $\bP^2$ and the Prym variety of a
naturally defined \'etale double cover of the discrminant curve of the
conic-bundle. Clemens and Griffiths used this isomorphism to give a
proof of the irrationality of a smooth cubic threefold and Beauville
later generalized the isomorphism to intermediate jacobians of
odd-dimensional quadric-bundles over $\bP^2$. We further generalize
the isomorphism to the primitive {\em cohomology} of a smooth cubic
hypersurface in $\bP^n$.
\end{abstract}
\maketitle

\vskip.5in
\begin{center}
{\sc Introduction}
\end{center}
\vskip.25in

Fano studied the variety of lines on a cubic hypersurface with a
finite number of singular points. The variety parametrizing linear
spaces of given dimension in a projective variety $X$ is now called a
Fano variety.  Subvarieties of a Fano variety can be defined using
various incidence relations. Such varieties are studied to help
understand the geometric properties of $X$ and for their own sake. For
instance, the proofs of the irrationality of a smooth cubic threefold
$X$ and of the Torelli theorem for $X$ by Clemens and Griffiths use
varieties of lines in the cubic.

Suppose that $X$ is a smooth cubic hypersurface in $\bP^4$ and let $F$
be the Fano variety of lines in $X$. By \cite{CG} (Lemma 7.7 page 312)
the variety $F$ is a smooth surface. Let us fix a general line $l$ in
$X$, corresponding to a general element of $F$, and let $D_l$ be the
variety of lines in $X$ incident to $l$. The blow up $X_l$ of $X$
along $l$ has the structure of a conic-bundle over $\bP^2$ and its
discriminant curve is a smooth plane quintic $Q_l$:
\[
\begin{array}{ccc}
 & & X_l \\
 & & \downarrow \\
Q_l & \subset & \bP^2 \: .
\end{array}
\]
The curve $D_l$
is an \'etale double cover of $Q_l$.

In a first proof of the irrationality of $X$, Clemens and Griffiths use the
canonical isomorphism between the Albanese variety of $F$ and the
intermediate jacobian of $X$ (see \cite{CG}, Theorem 11.19 page 334). In a
second proof they use the canonical isomorphism (due to Mumford, see
\cite{CG}, Appendix C) between the intermediate jacobian of $X$ and the
Prym variety of the (\'etale) double cover $D_l \rightarrow Q_l$. More
generally, Mumford's result says that this isomorphism holds for a
conic bundle $X$ over $\bP^2$ with discriminant curve $Q_l$ and double
cover $D_l$ parametrizing the components of the singular conics
parametrized by $Q_l$. Beauville generalized this isomorphism to the
case where $X$ is an odd-dimensional quadric bundle over $\bP^2$ with
discriminat curve $Q_l$ and double cover $D_l$ parametrizing the
rulings of the quadrics parametrized by $Q_l$ (see \cite{B2}).

In this paper we ``generalize'' the isomorphism between the intermediate
jacobian of $X$ and the Prym variety of $D_l \rightarrow Q_l$ to the
cohomology of higher dimensional cubic hypersurfaces. On the way we also
obtain some results about the Fano variety $\cP$ of planes in $X$.

A principally polarized abelian variety $A$ is the Prym variety of a double
cover of curves $\pi : \tC \rightarrow C$ if there is an exact sequence
\[
0 \lra \pi^* JC \lra J\tC \lra A \lra 0
\]
and, under the transpose of $J\tC \surj A$, the principal polarization
of $J\tC$ pulls back to twice the principal polarization of $A$. The
generalization that we have in mind would say that a polarized Hodge
structure $H$ is the Prym Hodge structure of two polarized Hodge
structures $H_1 \subset H_2$ if there is an involution $i : H_2
\rightarrow H_2$ and a surjective morphism of Hodge structures $\psi :
H_2 \surj H$ such that $i$ is a morphism of Hodge structures of type
$(0,0)$, the kernel of $\psi$ is the $i$-invariant part of $H_2$ which
is equal to $H_1$ and such that for any two $i$-anti-invariant
elements $a,b$ of $H_2$ we have $\psi(a).\psi(b)= - 2 a.b$ where ``.''
denotes the polarizations (see \cite{B2} page 334). In our case $H$
will be the primitive cohomology of a cubic hypersurface and $H_1$ and
$H_2$ will be the ``primitive'' cohomologies of (partial)
desingularizations of $Q_l$ and $D_l$.

From now on let $X$ be a smooth cubic hypersurface in $\bP^n$. For a
general line $l \subset X$, we define $X_l$ to be the blow up of $X$
along $l$. Then $X_l$ is a conic bundle over $\bP^{n-2}$ and we define
$Q_l$ to be its discriminant variety:
\[
\begin{array}{ccc}
 & & X_l \\
 & & \downarrow \\
Q_l & \subset & \bP^{n-2}\: .
\end{array}
\]
For $n \geq 5$ the variety $Q_l$ is singular. It parametrizes the
singular or higher-dimensional fibers of $X_l\rightarrow \bP^{n-2}$
and it can be thought of as the variety parametrizing planes in
$\bP^n$ which contain $l$ and, either are contained in $X$ or, whose
intersection with $X$ is a union of three (possibly equal) lines. We
define $D_l$ to be the variety of lines in $X$ incident to $l$. Then
$D_l$ admits a rational map of degree $2$ to $Q_l$ and the
varieties $D_l$ and $Q_l$ have dimension $n-3$. It is proved in
\cite{voisin}, page 590, that $D_l$ is smooth and its map to $Q_l$ is
a morphism for $n=5$ and $l$ general. For $n\geq 6$, the variety $D_l$
is always singular and for $n \geq 8$, the rational map $D_l
\rightarrow Q_l$ is never a morphism. We define a natural
desingularization $S_l$ of $D_l$ such that the rational map
$D_l\rightarrow Q_l$ lifts to a morphism $S_l\rightarrow Q_l$. The
morphism is not finite however. So we define natural blow-ups $S_l'$
and $Q_l'$ of $S_l$ and $Q_l$ such that the morphism $S_l\rightarrow
Q_l$ lifts to a double cover $S_l'\rightarrow Q_l'$. The varieties
$S_l$ and $S_l'$ naturally parametrize lines in blow-ups of $X_l$ so
that we have Abel-Jacobi maps $\psi : H^{n-3}(S_l,{\Bbb Z}) \lra
H^{n-1}(X,{\Bbb Z})$ and $\psi' : H^{n-3}(S_l',{\Bbb Z}) \lra
H^{n-1}(X,{\Bbb Z})$. Our main results are as follows
\begin{lemma0}
The Abel-Jacobi maps
\[ \psi : H^{n-3}(S_l,{\Bbb Z}) \lra H^{n-1}(X,{\Bbb Z}) \]
and
\[ \psi' : H^{n-3}(S_l',{\Bbb Z}) \lra H^{n-1}(X,{\Bbb Z}) \]
are onto.
\end{lemma0}
The involution $i_l : S_l'\rightarrow S_l'$ of the double cover
$S_l'\rightarrow Q_l'$ induces an involution $i : H^{n-3}(S_l',{\Bbb
Z})\rightarrow H^{n-3}(S_l',{\Bbb Z})$ whose invariant subgroup is
$H^{n-3}(Q_l',{\Bbb Z})$. However, the Prym construction only works
for ``primitive'' cohomologies (see Definition \ref{defprim}
below). Denote the primitive part of each cohomology group $H$ by
$H^0$. We need to show that for any two $i$-anti-invariant elements
$a,b$ of $H^{n-3}(S_l',{\Bbb Z})^0$, we have $\psi'(a).\psi'(b) = -2
a.b$. This follows from (see \ref{thmaibab})
\begin{theorem0}
Let $a$ and $b$ be two elements of $H^{n-3}(S_l', \bZ)^0$. Then
\[
\psi'(a). \psi'(b) = a. i_l^* b - a . b \: .
\]
\end{theorem0}
We use this theorem to prove
\begin{theorem0}
The Abel-Jacobi map
\[ {\psi'}^0 : H^{n-3}(S_l',{\Bbb Z})^0 \lra H^{n-1}(X,{\Bbb Z})^0 \]
is onto with kernel equal to the image of $H^{n-3}(Q_l', \bZ)^0$ in
$H^{n-3}(S_l',{\Bbb Z})^0$.
\end{theorem0}
This finishes the Prym construction.

Let us now briefly discuss the Hodge conjectures. The general Hodge
conjecture $GHC(X,m,p)$ as stated in \cite{steenbrink} page 166 is the
following

\begin{description}
\item[$GHC(X,m,p)$] For every $\bQ$-Hodge substructure $V$ of $H^m(X,
\bQ)$ with level $\leq m - 2p$, there exists a subvariety $Z$ of $X$
of codimension $p$ such that $V$ is contained in the image of the
Gysin map $H^{m-2p}(\tZ, \bQ) \lra H^m(X, \bQ)$ where $\tZ$ is a
desingularization of $Z$.
\end{description}

It is proved in \cite{steenbrink}, Proposition 2.6, that $GHC(Y,m,1)$
holds for all uniruled smooth varieties $Y$ of dimension $m$. Our
Lemma $1$ gives a geometric proof of GHC(X,n-1,1) for a smooth cubic
hypersurface $X$ in $\bP^n$: the full cohomology $H^{n-1}(X,
\bZ)$ is supported on the subvariety $Z$ which is the union of all the
lines in $X$ incident to $l$.

We now describe our results in slightly greater detail. In Section $1$
we prove that, for $n \geq 6$ and $l$ general, the singular locus of
$D_l$ is $\{ l \} \subset D_l$. Also, for $n \geq 6$, the natural map
$D_l \rightarrow Q_l$ sending a line $l'$ to the plane spanned by $l$
and $l'$ is only a rational map. In Section $2$, we prove that the
variety $S_l$ parametrizing lines in the fibers of the conic bundle
$X_l\rightarrow \bP^{n-2}$ is a small desingularization of $D_l$ which
admits a {\em morphism} of generic degree $2$ to $Q_l$. We show that
$S_l$ can also be defined as a subvariety of the product of
Grassmannians of lines and planes in $\bP^n$. For the generalized Prym
construction we need a finite morphism of degree $2$ to $Q_l$ and the
morphism $S_l \rightarrow Q_l$ is not finite for $n \geq 8$. It fails
to be finite at the points of $Q_l$ parametrizing planes contained in
$X$ (and containing $l$). Let $T_l$ denote the variety parametrizing
planes in $X$ which contain $l$. Since $\bP^{n-2}$ parametrizes the planes
in $\bP^n$ which contain $l$, the variety $T_l$ is naturally a
subvariety of $\bP^{n-2}$ and in fact is contained in $Q_l$:
\[
\begin{array}{ccc}
 & & X_l \\
 & & \downarrow \\
T_l \subset Q_l & \subset & \bP^{n-2}\: .
\end{array}
\]
In Section $3$ we prove that for $l$ general $T_l$ is a smooth
complete intersection of the expected dimension $n-8$ in
$\bP^{n-2}$. For this we analyse the structure of the Fano variety
$\cP$ of planes in $X$. We prove that $\cP$ is always of the expected
dimension $3n-16$ and determine its singular locus. It is proved in
\cite{borcea1} Theorem 4.1 page 33 or \cite{debarremanivel}
Th\'eor\`eme 2.1 pages 2-3 that $\cP$ is connected for $n \geq 6$. We
prove that $\cP$ is irreducible for $n \geq 8$. In Section $4$ we blow
up $X_l \rightarrow \bP^{n-2}$ along $T_l$ and its inverse image in
$X_l$ to obtain $X_l' \rightarrow {\bP^{n-2}}'$. The discriminant
hypersurface of this conic-bundle is the blow-up $Q_l'$ of $Q_l$ along
$T_l$:
\[
\begin{array}{ccc}
 & & X_l' \\
 & & \downarrow \\
Q_l' & \subset & {\bP^{n-2}}'\: .
\end{array}
\]
The variety $S_l'$ is then defined as the variety of lines in the
fibers of the conic bundle $X_l' \rightarrow {\bP^{n-2}}'$. We prove
that the rational involution acting in the fibers of $S_l \rightarrow
Q_l$ lifts to a regular involution $i_l : S_l' \rightarrow S_l'$ and
the quotient of $S_l'$ by $i_l$ is $Q_l'$. We also prove that $S_l'$
is the blow up of $S_l$ along the inverse image of $T_l$ and the
ramification locus $R_l'$ of $S_l'\rightarrow Q_l'$ is smooth of
codimension $2$ and is an ordinary double locus for $Q_l'$. In Section
$5$ we prove Lemma 1, Theorem 2 and Theorem 3. We also prove some
results about the rational cohomology ring of $S_l$: we prove that,
except in the middle degree, this rational cohomology ring is
generated by $H^2(S_l, \bQ)$ which, for $n \geq 6$, is generated by
the inverse images $h$ and $\sigma_1$ of the hyperplane classes of
$Q_l$ and $D_l$ (the hyperplane class of $D_l$ is the restriction of
the hyperplane class of the Grassmannian of lines in $\bP^n$). For
$n=5$, the space $H^2(S_l, \bQ)$ is the direct sum of its primitive
part and $\bQ h\oplus\bQ\sigma_1$.

{\sc Acknowledgments}

I am indebted to Joe Harris, Brendan Hassett, Barry Mazur and Robert
Varley for many stimulating discussions.

\vskip.5in
\begin{center}
{\sc Notation and Conventions}
\end{center}
\vskip.25in

The symbol $n$ will always denote an integer $\geq 5$.

For all positive integers $k$ and $l$, we denote by $G(k,l)$ the
Grassmannian of $k$-dimensional vector spaces in $\bC^l$. For any vector 
space or vector bundle $W$, we denote by $\bP (W)$ the projective space 
of lines in (the fibers of) $W$ with its usual scheme structure.

For all cohomology vector spaces $H^i(Y, . )$ of a variety $Y$, we will
denote by $h^i(Y, . )$ the dimension of $H^i(Y, .)$. For a point $y
\in Y$, we denote by $T_yY$ the Zariski tangent space to $Y$ at $y$.
If we are given an embedding $Y \subset \bP^m  = (\bC^{m+1} \setminus
\{ 0 \})/ \bC^*$, we denote by $T'_yY$ the inverse image of $T_yY$ by
the map $\bC^{m+1} \surj \bC^{m+1} / \bC v = T_y \bP^m $ where $v$
is a nonzero vector in $\bC^{m+1}$ mapping to $y$. We call $\bP(T'_yY)$ the
projective tangent space to $Y$ at $y$.

For any subsets or subschemes $Y_1, ..., Y_m$ of a projective space
$\bP^d$ (resp. affine space $\bC^d$), we denote by $\langle Y_1, ...,
Y_m \rangle$ the smallest linear subspace of $\bP^d$ (resp. $\bC^d$)
containing $Y_1, ..., Y_m$.

For a subscheme $Y_1$ of a scheme $Y_2$, we denote by $N_{Y_1/Y_2}$ the 
normal sheaf to $Y_1$ in $Y_2$.

For a global section $s$ of a sheaf $\cF$ on a scheme $Y$, we denote by
$Z(s)$ the scheme of zeros of $s$ in $Y$.

\section{The variety $D_l$ of lines incident to $l$}
\label{secDl}

For a smooth cubic hypersurface $X \subset \bP^n$ of equation $G$, we let
$\delta : \bP^n \lra (\bP^n)^*$ be the dual morphism of $X$. In terms of a
system of projective coordinates $\{x_0, ..., x_n \}$ on $\bP^n$, the
morphism $\delta$ is given by
\[
\delta (x_0, ..., x_n) = \left( \del_0 G(x_0, ..., x_n), ..., \del_n G(x_0,
..., x_n) \right)
\]
where $\del_i = \frac{\del}{\del x_i}$.

Let $l \subset X$ be a line. Following \cite{CG} (page 307 Definition
6.6, Lemma 6.7 and page 310 Proposition 6.19), we make the
definition:
\begin{definition}
\begin{enumerate}
\item The line $l$ is of first type if the normal bundle to $l$ in $X$ is
isomorphic to $\cO_l^{\oplus 2} \oplus \cO_l(1)^{\oplus
(n-4)}$. Equivalently, the intersection $\bT_l$ of the projective tangent
spaces to $X$ along $l$ is a linear subspace of $\bP^n$ of dimension
$n-3$. Equivalently, the dual morphism $\delta$ maps $l$ isomorphically
onto a conic in $(\bP^n)^*$, i.e., the restriction map $\langle \del_0 G,
..., \del_n G \rangle \lra H^0(l, \cO_l(2))$ is onto where $\langle \del_0
G, ..., \del_n G \rangle$ is the span of $\del_0 G, ..., \del_n G$ in
$H^0(\bP^n, \cO_{\bP^n}(2))$.
\item The line $l$ is of second type if the normal bundle to $l$ in $X$ is
isomorphic to $\cO_l(-1) \oplus \cO_l(1)^{\oplus (n-3)}$. Equivalently, the
space $\bT_l$ is a linear subspace of $\bP^n$ of dimension
$n-2$. Equivalently, the dual morphism $\delta$ has degree $2$ on $l$ and
maps $l$ onto a line in $(\bP^n)^*$, i.e., the restriction map $\langle
\del_0 G, ..., \del_n G \rangle \lra H^0(l, \cO_l(2))$ has rank $2$.
\end{enumerate}
\end{definition}

By \cite{CG} (Lemma 7.7 page 312), the variety $F$ of lines in $X$ is
smooth of dimension $2(n-3)$. An easy dimension count shows that the
dimension of $D_l$ is at least $n-3$ for any $l \in F$. Suppose that
$l$ is of first type. We have

\begin{lemma}
Let $l' \in D_l$ be distinct from $l$. If $l'$ is of first type or if $l'$
is of second type and $l$ is {\em not} contained in $\bT_{l'}$, then the
dimension of $T_{l'}D_l$ is $n-3$ (i.e., $D_l$ is smooth of dimension $n-3$
at $l'$). If $l'$ is of second type and $l$ is contained in $\bT_{l'}$, then
the dimension of $T_{l'}D_l$ is $n-2$.
\label{Dlsmooth1}
\end{lemma}
{\em Proof :} The variety $D_l$ is the intersection of $F$ with the variety
$G_l$ parametrizing all lines in $\bP^n$ which are incident to
$l$. Therefore $T_{l'}D_l = T_{l'}G_l \cap T_{l'}F \subset
T_{l'}G(2,n+1)$. Let $V$ and $V'$ be the vector spaces in ${\Bbb C}^{n+1}$
whose projectivizations are respectively $l$ and $l'$. Then $T_{l'}G_l$ can
be identified with the subvector space of $T_{l'}G(2,n+1) = Hom(V', {\Bbb
C}^{n+1} / V')$ consisting of those homomorphisms $f$ such that $f(V \cap
V') \subset (V+V') / V'$ (see e.g. \cite{harris2}, Ex. 16.4 pages
202-203). It follows that the set of homomorphisms $f$ such that $f(V \cap
V') = 0 $ is a codimension one subspace of $T_{l'}G_l$ and therefore its
intersection $H$ with $T_{l'}D_l$ has codimension one or less in
$T_{l'}D_l$.

The space $T_{l'}F$ can be identified with the subvector space of
$T_{l'}G(2,n+1) = Hom(V', {\Bbb C}^{n+1} / V')$ consisting of those
homomorphisms $f$ such that for any vector $v \in V' \setminus \{ 0 \}$
mapping to a point $p \in l'$, we have $f(v) \in T'_pX / V'$ (see
\cite{harris2} Ex. 16.21 and 16.23 pages 209-210). If $f : V' \lra {\Bbb
C}^{n+1} / V'$ verifies $f(V \cap V') =0$, then $f(V') = {\Bbb C} f(v)$ for
$v$ a general vector in $V'$. Hence, if $f \in H$, then $f(V') \subset
\bigcap_{p \in l'} T_p'X / V'$.

If $l'$ is of first type, then $\bigcap_{p \in l'} T_p'X$ has dimension
$n-2$, hence $\bigcap_{p \in l'} T_p'X / V'$ has dimension $n-4$. So $H$
has dimension $n-4$ and, since $H$ has codimension one or less in
$T_{l'}D_l$, we deduce that $T_{l'}D_l$ has dimension at most $n-3$ hence
it has dimension equal to $n-3$ (since $D_l$ has dimension $\geq n-3$).

If $l'$ is of second type, then the tangent space $T_{l'}F$ can be
identified with $Hom(V', \bigcap_{p \in l'} T_p'X / V')$ (because, for
instance, the latter is contained in $T_{l'}F$ and the two spaces have the
same dimension). If $V$ is not contained in $\bigcap_{p \in l'} T_p'X$,
then $f(V \cap V') \subset (V+V') / V'$ for $f \in Hom(V', \bigcap_{p \in
l'} T_p'X / V')$ implies $f(V \cap V') =0$. So $T_{l'}D_l = T_{l'}F \cap
T_{l'}G_l$ has dimension equal to the dimension of $\bigcap_{p \in l'}
T_p'X / V'$ which is $n-3$. So in this case $D_l$ is smooth at $l'$.  If $V
\subset \bigcap_{p \in l'} T_p'X / V'$, then the requirement $f(V \cap V')
\subset (V+V') / V'$ imposes $n-4$ conditions on $f$ and the dimension of
$T_{l'}D_l$ is $n-2$.  \hfill $\qed$ \vskip10pt

Since $\bT_l$ has dimension $n-3$, we see that, as soon as $n \geq 5$, we
have $l \in D_l$. We have the following

\begin{lemma}
If $n \geq 6$, then $D_l$ is singular at $l$. If $n = 5$, then $D_l$
is smooth at $l$ if $X$ does not have contact multiplicity $3$ along
$l$ with the plane $\bT_l$ and if there is no line $l'$ of second type
in $\bT_l$.
\end{lemma}
{\em Proof :} The case $n=5$ is Lemme 1 on page 590 of
\cite{voisin}. Suppose $n \geq 6$. For $l$ general, consider a plane
section of $X$ of the form $l+l'+l''$ such that $l \cap l'$ and $l \cap
l''$ are general points on $l$. The set of lines through $l \cap l'$ is a
divisor in $D_l$ and meets the set of lines through $l \cap l''$ only at $l
\in D_l$. So we have two divisors in $D_l$ which meet only at a point and
$D_l$ has dimension $\geq 3$. Therefore $D_l$ is not smooth at $l$ for $l$
general and hence for all $l$.  \hfill $\qed$ \vskip10pt

We now prove an existence result, namely,
\begin{lemma}
The set of lines $l \in F$ such that $l$ is contained in $\bT_{l'} $ for
some line $l' \in F$ of second type is a proper closed subset of $F$. In
other words (by Lemma \ref{Dlsmooth1}), for $l \in F$ general, the variety
$D_l \setminus \{ l \}$ is smooth of dimension $n-3$.
\label{lemTl'}
\end{lemma}

{\em Proof :} Since the dimension of $F$ is $2(n-3)$ and the dimension
of the variety $F_0 \subset F$ parametrizing lines of second type is
$n-3$ (see \cite{CG} page 311 Corollary 7.6), if the lemma fails, then
for any line $l' \in F_0$, the dimension of the family of lines
in $X \cap \bT_{l'}$ which intersect $l'$ is at least $n-3$.

The variety $\bT_{l'}$ is a linear subspace of codimension $2$ of
$\bP^n$. Any plane in $\bT_{l'}$ which contains $l'$ is tangent to $X$
along $l'$. The intersection of a general such plane $P$ with $X$ is the
union of $l'$ and a line $l$, the line $l'$ having multiplicity $2$ (or $3$
if $l=l'$) in the intersection cycle $[P \cap X]$. Conversely, any line $l$
in $X \cap \bT_{l'}$ which intersects $l'$ is contained in such a
plane. The family of planes in $\bT_{l'}$ which contain $l'$ has dimension
$n-4$. Therefore, if the family of lines $l$ in $X \cap \bT_{l'}$ which
intersect $l'$ has dimension $\geq n-3$, then for each such line $l \neq
l'$, the plane $\langle l, l' \rangle$ contains a positive-dimensional
family of lines in $X \cap \bT_{l'}$ and hence $\langle l, l' \rangle$ is
contained in $X \cap \bT_{l'}$. Therefore $X \cap \bT_{l'}$ is a cone over
a cubic hypersurface in $\bT_{l'} / l'$ and, for each plane $P \subset X
\cap \bT_{l'}$ which contains $l'$, there is a hyperplane in $\bT_{l'}$
tangent to $X \cap \bT_{l'}$ along $P$. Therefore $\bT_P \stackrel{def}{=}
\cap_{p \in P} \bP T'_p X$ has codimension $3$ in $\bP^n$. Hence the
restriction of the dual morphism of $X$ to $P$ is a morphism of degree $4$
from $P$ onto a plane in $(\bP^n)^*$. It follows from \cite{CG} Lemma 5.15
page 304 that all such planes are contained in a proper closed subset of
$X$. Therefore a general line $l \in F$ is not contained in such a plane
and hence not in $\bT_{l'}$. Contradiction. \hfill $\qed$

\section{Desingularizing $D_l$}
\label{secSl}

Let $X_l$ and $\bP^n_l$ be the blow ups of $X$ and $\bP^n$
respectively along $l$. Then the projection from $l$ gives a
projective bundle structure on $\bP^n_l$ and a conic bundle structure
on $X_l$ (i.e., a general fiber of $\pi_X : X_l \lra \bP^{n-2}$ is a
conic in the corresponding fiber of $\pi : \bP^n_l \lra \bP^{n-2}$):
\[ \begin{array}{crl}
X_l & \inj & \bP^n_l \\
 & \pi_X \searrow & \downarrow \pi \\
 & & \bP^{n-2}
\end{array}
\]
Let $E$ be the locally free sheaf $\cO_{{\Bbb P}^{n-2}}(-1) \oplus
\cO_{\bP^{n-2}}^{\oplus 2}$. Then it is easily seen (as in
e.g. \cite{hartshorne} page 374 Example 2.11.4) that $\pi : \bP^n_l \lra
\bP^{n-2}$ is isomorphic to the projective bundle $\bP(E) \lra \bP^{n-2}$.
The variety $X_l \subset \bP^n_l$ is the divisor of zeros of a section $s$
of $\cO_{\bP E}(2) \otimes \pi^* \cO_{\bP^{n-2}}(m)$ for some integer $m$
because the general fibers of $\pi_X : X_l \lra \bP^{n-2}$ are smooth
conics in the fibers of $\pi$. Since $\pi_* (\cO_{\bP E}(2) \otimes \pi^*
\cO_{{\Bbb P}^{n-2}}(m)) \cong Sym^2E^* \otimes \cO_{\bP^{n-2}}(m)$, the
section $s$ defines a (``symmetric'') morphism of vector bundles $\phi : E
\lra E^* \otimes \cO_{\bP^{n-2}}(m)$. The degeneracy locus $Q_l \subset
{\Bbb P}^{n-2}$ of this morphism is the locus over which the fibers of
$\pi_X$ are singular conics (or have dimension $\geq 2$). By, for instance,
intersecting $Q_l$ with a general line, we see that $Q_l$ is a quintic
hypersurface (see \cite{segre} pages 3-5). Therefore $m=1$. Let $S_l$ be
the variety parametrizing lines in the fibers of $\pi_X$. We have a
morphism $S_l \lra D_l$ defined by sending a line in a fiber of $\pi$ to
its image in $\bP^n$. Let $E_1 \subset X_l$ be the exceptional divisor of
$\epsilon_1 : X_l \rightarrow X$ and let $P_1 \subset S_l$ be the variety
parametrizing lines which lie in $E_1$. Then the morphism $S_l \lra D_l$
induces an isomorphism $S_l \setminus P_1 \cong D_l \setminus \{ l \} $. We
have

\begin{lemma}
Suppose that $l$ is of first type and $D_l \setminus \{ l \}$ is
smooth. Then $S_l$ is smooth and irreducible and admits a {\em morphism} of
generic degree $2$ onto $Q_l$. The variety $S_l$ can also be defined as the
closure of the subvariety of $G(2,n+1) \times G(3,n+1)$ parametrizing pairs
$(l',L')$ such that $l' \in D_l \setminus \{ l \}$ and $L' = \langle l,l'
\rangle$.
\label{lemSlsm}
\end{lemma}
{\em Proof :} The morphism $S_l \lra Q_l$ is defined by sending a line
in a fiber of $\pi$ to its image in $\bP^{n-2}$. It is of generic
degree $2$ because the rational map $D_l \lra Q_l$ is of generic
degree $2$. The variety $S_l$ is irreducible because $Q_l$ is irreducible
and $S_l \lra Q_l$ is not split (intersect $Q_l$ with a general plane and
use \cite{B5} ).

For $l' \in S_l \setminus P_1$, the variety $S_l$ is smooth at $l'$
since $S_l \setminus P_1 \cong D_l \setminus \{ l \} $.

For $l' \in P_1$ we determine the Zariski tangent space to $S_l$ at
$l'$. Since $l'$ maps to a point in $\bP^{n-2}$, it corresponds to a plane
$L'$ in ${\Bbb P}^n$ which contains $l$. Since $l'$ is also contained in
$E_1$, it maps onto $l$ in $\bP^n$ under the blow up morphism $\bP^n_l \lra
\bP^n$ and $L'$ is tangent to $X$ along $l$. So we easily see that we can
identify $S_l$ with the closure of the subvariety of the product of the
Grassmannians $G(2,n+1) \times G(3,n+1)$ parametrizing pairs $(l',L')$ such
that $l' \in D_l \setminus \{ l \}$ and $L' = \langle l,l' \rangle$.

Let $W'$ and $V$ be the vector spaces in $\bC^{n+1}$ whose
projectivizations are $L'$ and $l$ respectively. The tangent space to
$G(2,n+1) \times G(3,n+1)$ at $(l,L')$ can be canonically identified with $
Hom(V, {\Bbb C}^{n+1}/V) \oplus Hom(W', {\Bbb C}^{n+1}/W') $. As in
\cite{harris2} Ex. 16.3 pages 202-203 and Ex. 16.21, 16.23 pages 209-210,
one can see that the tangent space to $S_l$ at $(l,L')$ can be identified
with the set of pairs of homomorphisms $(f,g)$ such that for every nonzero
vector $v \in V$ mapping to a point $p$ of $l$, we have $f(v) \in T'_pX /
V$, $g(V) = 0$, $g |_V = f(\hbox{mod} W')$ and $g(W') \subset \cap_{p \in
l} T'_pX /W'$ (this last condition expresses the fact that the deformation
of $L'$ contains a deformation of $l$ which is contained in $X$, hence the
deformation of $L'$ is tangent to $X$ along $l$, i.e., is contained in
$\bT_l$). Equivalently, $g(V) =0$, $f(V) \subset W' /V$ and $g(W') \subset
\cap_{p \in l} T'_pX /W'$. Assuming $l$ of first type, we see that the
space of such pairs of homomorphisms has dimension $n-3$.  \hfill $\qed$

\section{The planes in $X$}
\label{secplanes}

Let $\cP$ be the variety parametrizing planes in $X$. For $P
\in \cP$, we say that $\delta$ has rank $r_P$ on $P$ if the span of
$\delta(P)$ has dimension $r_P$. Since $\delta$ is defined by
quadrics, we have $r_P \leq 5$. Since $X$ is smooth, we have $r_P \geq
2$. Consider the commutative diagram
\[
\begin{array}{rcl}
 & & \bP^5 \\
 & v \nearrow & \downarrow p \\
P & \stackrel{\delta_P}{\lra} & \bP^{r_P} \subset (\bP^n)^*
\end{array}
\]
where $v$ is the Veronese map, $\delta_P$ is the restriction of
$\delta$ to $P$ and $p$ is the projection from a linear space $L
\subset \bP^5$ of dimension $4 - r_P$ (with the convention that the
empty set has dimension $-1$).

Note that $L$ does not intersect $v(P)$ because $\delta$ is a morphism.

Let $\cP_r$ be the subvariety of $\cP$ parametrizing planes $P$ for which 
$r_P \leq r$. In this section we will prove a few facts about $\cP$ and 
$\cP_r$ which we will need later. We begin with

\begin{lemma}
Let $T \stackrel{def}{=} \cup_{l \subset P} \langle v(l) \rangle
\subset \bP^5$ be the secant variety of $v(P)$. Then there is a
bijective morphism from $T \cap L$ to the parameter space of the
family of lines of second type in $P$ and $T \cap L$ contains no
positive-dimensional linear spaces. In particular,
\begin{enumerate}
\item if $r_P = 5$, then $P$ contains no lines of second type,
\item if $r_P=4$, then $P$ contains at most one line of second type and
this happens exactly when $L$ (which is a point in this case) is in $T$,
\item if $r_P = 3$, then $P$ contains one, two or three
distinct lines of second type,
\item if $r_P=2$, then $P$ contains exactly a one-parameter family of lines
of second type whose parameter space is the bijective image of an
irreducible and reduced plane cubic.
\end{enumerate}
\label{linesinP}
\end{lemma}
{\em Proof :} A line $l \subset P$ is of second type if and only if
$\delta_P(l) \subset \bP^{r_P}$ is a line, i.e., if and only if the span
$\langle v(l) \rangle \cong \bP^2$ of the smooth conic $v(l)$ intersects
$L$. Consider the universal line $f_1: L_P \rightarrow P^*$ and its
embedding $L_P \inj V_P$ where $f_2:V_P \rightarrow P^*$ is the
projectivization of the bundle $f_* \cO_{L_P}(2)^*$. Then $T$ is the image
of $V_P$ in $\bP^5$ by a morphism, say $g$, which is an isomorphism on the
complement of $L_P$ and contracts $L_P$ onto $v(P)$. Since $L \cap v(P) =
\emptyset$, the morphism $g |_{g^{-1}(T \cap L)}$ is an isomorphism, say
$g'$. The morphism from $T \cap L$ onto the parameter space of the family
of lines of second type in $P$ is the composition of ${g'}^{-1}$ with
$f_2$. This morphism is bijective because (since $L \cap v(P) = \emptyset$)
the space $L$ intersects any $\langle v(l) \rangle$ in at most one point
and any two planes $\langle v(l_1) \rangle$ and $\langle v(l_2) \rangle$
intersect in exactly one point which is $v(l_1 \cap l_2) \in v(P)$.

To show that $T \cap L$ contains no positive-dimensional linear spaces,
recall that $T$ is the image of the Segre embedding of $P \times P$ in
$\bP^8 = \bP \left( H^0(P, \cO_P (1))^{\otimes 2} \right)^*$ by the
projection from $\bP \left( \Lambda^2 H^0(P, \cO_P (1)) \right)^*$. Let
$R_1$ be the ruling of $T$ by planes which are images of the fibers of the
two projections of $P \times P$ onto $P$. Let $R_2$ be the ruling of $T$ by
planes of the form $\langle v(l) \rangle$ for some line $l \subset P$. Then
a simple computation (determining all the pencils of conics which consist
entirely of singular conics) shows that every linear subspace contained in
$T$ is contained either in an element of $R_1$ or an element of
$R_2$. Therefore, if $L \cap T$ contains a linear space $m$, then either $m
\subset \langle v(l) \rangle$ for some line $l \subset P$ or $m \subset L'$
for some element $L'$ of $R_1$. In the first case, the space $m$ is a point
because otherwise it would intersect $v(P)$. In the second case, the space
$m$ is either a point or a line because any element of $R_1$ contains
exactly one point of $v(P)$. It is easily seen that there is an element
$s_0 \in H^0(P,\cO_P(1))$ such that $L'$ parametrizes the hyperplanes in
$|\cO_P(2)|$ containing all the conics of the form $Z(s.s_0)$ for some $s
\in H^0(P,\cO_P(1))$. If $m \subset L'$ is a line, then it is easily seen
that the codimension, in $\langle \del_0 G, ..., \del_n G \rangle |_P$, of
the set of elements of the form $s.s_0$ is one. Restricting to $Z(s_0)$, we
see that the dimension of $\langle \del_0 G, ..., \del_n G \rangle
|_{Z(s_0)}$ is $1$ which is impossible since then $X$ would have a singular
point on $Z(s_0)$. Therefore $m$ is always a point if non-empty.  \hfill
$\qed$

\begin{proposition}
The space of infinitesimal deformations of $P$ in $X$
has dimension $3n-15$ if $r_P = 2$. In particular, if $n=5$, then $X$
contains at most a finite number of planes.
\label{rP2inf}
\end{proposition}
{\em Proof :} The intersection $\bT_P$ of the projective tangent spaces to
$X$ along $P$ has dimension $n-3$. It follows that we have an exact
sequence
\[
0 \lra \cO_P(1)^{n-5} \lra N_{P/X} \lra V_2 \lra 0
\]
where $V_2$ is a locally free sheaf of rank $2$. We need to show that
$h^0(P,V_2) = 0$. Suppose that there is a nonzero section $u
\in H^0(P,V_2)$. We will first show that the restriction of $u$ to any
line $l$ in $P$ is nonzero. This will follow if we show that the
restriction map $H^0(P,V_2) \lra H^0(l, V_2|_l)$ is injective, i.e.,
$h^0(P, V_2(-1)) = 0$. Consider therefore the exact sequence of normal
sheaves
\[ 0 \lra N_{P/X} \lra N_{P/\bP^n} \lra N_{X/\bP^n} |_P \lra 0 \]
After tensoring by ${\cal
O}_P(-1)$ we obtain the exact sequence
\[ 0 \lra N_{P/X}(-1) \lra {{\cal O}_P}^{\oplus (n-2)} \lra {\cal O}_P(2)
\lra 0.
\]
We can choose our system of coordinates (on $\bP^n$) in such a way that
$x_3 = ... = x_n = 0$ are the equations for $P$ and the map ${{\cal
O}_P}^{\oplus (n-2)} \lra {\cal O}_P(2)$ in the sequence above is given by
multiplication by $\del_3 G |_P, ... \del_n G |_P$. So we see that, since
$r_P = 2$, the map on global sections $H^0({{\cal O}_P}^{\oplus (n-2)})
\lra H^0({\cal O}_P(2))$ has rank $3$. Therefore $h^0(P,N_{P/X}(-1)) = n-5$
and $h^0(P, V_2(-1)) = 0$.

By Lemma \ref{linesinP}, the plane $P$ contains lines of first
type. For any line $l \subset P$ which is of first type, it is easily
seen that $V_2 |_l \cong \cO_l^{\oplus 2}$. Hence $u$ has no zeros on
$l$. It follows that $Z(u)$ is finite.

We compute the total Chern class of $V_2$ as
\[ c(V_2) = \frac{c(N_{P/X})}{(1+ \zeta)^{n-5}} = 1+3 \zeta^2 \]
where $\zeta = c_1(\cO_P(1))$. Therefore $Z(u)$ is a finite subscheme of
length $3$ of $P$. Let $l_u$ be a line in $P$ such that $l_u \cap
Z(u)$ has length $\geq 2$. Then, by what we saw above, $l_u$ is of
second type. It is easily seen that $V_2 |_{l_u} \cong \cO_{l_u}(-1)
\oplus \cO_{l_u}(1)$. Restricting $u$ to $l_u$, we see that
$Z(u|_{l_u})= l_u \cap Z(u)$ has length $1$ which is a
contradiction. So $h^0(P,V_2)=0$ and $h^0(P,N_{P/X}) = 3n-15$.
\hfill
$\qed$

\vskip20pt

The next result we will need is

\begin{lemma}
The dimension of $\cP_2$ is at most $Min(n-4, 5)$.
\label{rP2}
\end{lemma}
{\em Proof :} The proof of the part $dim(\cP_2) \leq n-4$ is
similar to the proof of Corollary 7.6 on page 311 of \cite{CG}.

To prove $dim(\cP_2) \leq 5$, we may suppose that $n \geq 10$. Let $P$ be
an element of $\cP_2$. We will show that the space of infinitesimal
deformations of $P$ for which the rank of $\delta$ does not increase has
dimension at most $5$. Let $x_0, x_1, x_2$ be coordinates on $P$, let $x_0,
x_1, x_2, x_3, ..., x_{n-3}$ be coordinates on $\bT_P$ and $x_0, ...,
x_{n-3}, x_{n-2}, x_{n-1}, x_n$ coordinates on $\bP^n$. Then the conditions
$P \subset X$ and $\bT_P$ is tangent to $X$ along $P$ can be written
\[
\partial_i \partial_j \partial_k G = 0
\]
for all $i,j \in \{ 0,1,2 \}$, $k \in \{ 0, ..., n-3 \}$, where $G$
is, as before, an equation for $X$ and $\del_i= \frac{\del}{\del
x_i}$. We need to determine the infinitesimal deformations of $P$ for
which there is an infinitesimal deformation of $\bT_P$ which is
tangent to $X$ along the deformation of $P$. The infinitesimal
deformations of $P$ in $\bP^n$ are parametrized by
\[
Hom_{\bC} \left( \langle \partial_0, \partial_1, \partial_2 \rangle,
\frac{\bC^{n+1}}{\langle \partial_0, \partial_1, \partial_2 \rangle}
\right) \cong Hom_{\bC} \left( \langle \partial_0, \partial_1,
\partial_2 \rangle, \langle \partial_3, ..., \partial_n \rangle \right)
\]
and those of $\bT_P$ in $\bP^n$ are parametrized by
\[
Hom_{\bC} \left( \langle \partial_0, ..., \partial_{n-3} \rangle,
\frac{\bC^{n+1}}{\langle \partial_0, ..., \partial_{n-3} \rangle} \right) 
\cong
Hom_{\bC} \left( \langle \partial_0, ..., \partial_{n-3} \rangle, \langle
\partial_{n-2}, \del_{n-1}, \partial_n \rangle \right)
\]
where we also denote by $\partial_i$ the vector in $\bC^{n+1}$
corresponding to the differential operator $\partial_i$. We need to
determine the homomorphisms $\{ \partial_i \mapsto \partial'_i \in
\langle \partial_3, ..., \partial_n \rangle : i \in \{ 0, 1, 2 \} \}$
for which there is a homomorphism $\{ \partial_i \mapsto \partial''_i
\in \langle \del_{n-2}, \del_{n-1}, \del_n \rangle : i \in \{ 0, ...,
n-3 \} \}$ such that the following conditions hold.

\begin{enumerate}
\item The vector $\partial''_i$ is the projection of $\partial_i'$ to
$\langle \del_{n-2}, \del_{n-1}, \del_n \rangle$ for $i \in \{0,1,2
\}$. This expresses the condition that the infinitesimal deformation
of $\bT_P$ contains the infinitesimal deformation of $P$.

\item For all $i,j \in \{0,1,2 \}$ and $k \in \{0, ..., n-3 \}$,
\[
\left( \partial_i + \epsilon \partial'_i \right) \left( \partial_j +
\epsilon \partial'_j \right) \left( \partial_k + \epsilon \partial''_k
\right) G = 0
\]
where, as usual, $\epsilon$ has square $0$. Here we are
``differentiating'' the relations $\partial_i \partial_j \partial_k G
=0$. Developing, we obtain
\[
\left( \partial_i \partial_j \partial_k'' + \partial_i \partial_j' 
\partial_k + \partial_i' \partial_j
\partial_k \right) G = 0 \: \: .
\]
\end{enumerate}
\vskip10pt
Writing $\partial_i' = \sum_{j=3}^n a_{ij} \partial_j$ and $\partial_i'' =
\sum_{j=n-2}^n b_{ij} \partial_j$, the above conditions can be written as
\begin{enumerate}
\item For all $i \in \{0,1,2 \}$ and $j \in \{n-2, n-1, n \}$,
\[
a_{ij} = b_{ij} \: \: .
\]

\item For all $i,j \in \{0,1,2 \}$ and $k \in \{0, ..., n-3 \}$,
\[
\sum_{l=n-2}^n b_{kl} \partial_i \partial_j \partial_l G + \sum_{l=3}^n 
a_{jl}
\partial_i \partial_l \partial_k  G + \sum_{l=3}^n a_{il} \partial_l 
\partial_j \partial_k G =
0 \: \: .
\]
\end{enumerate}
Incorporating the first set of conditions in the second and using the
relations $\del_i \del_j \del_k G = 0$ for $i,j \in \{0,1,2 \}, k \in \{ 0,
... n-3 \}$, we divide our
conditions into two different sets of conditions as follows. We are
looking for matrices $(a_{il})_{0 \leq i \leq 2, 3 \leq l \leq n}$ for
which there is a matrix $(b_{kl})_{3 \leq k \leq n-3, n-2 \leq l \leq
n}$ such that, for all $i,j,k \in \{ 0,1,2 \}$,
\[
\sum_{l=n-2}^n \left( a_{kl} \partial_i \partial_j \partial_l + a_{jl}
\partial_i \partial_l \partial_k + a_{il} \partial_l \partial_j
\partial_k \right) G = 0
\]
and, for all $i,j \in \{ 0,1,2 \}$, $k \in \{ 3, ..., n-3 \}$,
\[
\sum_{l=n-2}^n b_{kl} \partial_i \partial_j \partial_l G +
\sum_{l=3}^n \left( a_{jl}
\partial_i \partial_l \partial_k + a_{il} \partial_l \partial_j
\partial_k \right) G = 0 \:\: .
\]

Consider the matrix whose columns are indexed by the $a_{lm}, b_{su}$ ($0
\leq l \leq 2$, $3 \leq m \leq n$, $3 \leq s \leq n-3$, $n-2 \leq u \leq
n$), whose rows are indexed by {\em unordered} triples $(i,j,k)$ with $i,j
\in \{ 0,1,2 \}$, $k \in \{ 0, ..., n-3 \}$ and whose entries are the
$\del_i \del_j \del_m G$, $\del_i \del_m \del_k G$, $\del_m \del_j \del_k
G$ or $\del_i \del_j \del_u G$. The entry in the column of $a_{lm}$ and the
row of $(i,j,k)$ is nonzero only if $l=i,j$ or $k$. We can, and will,
suppose that $l=i$. Here is the list of possibly nonzero such entries.
\[
3 \leq m \leq n \:\: , \:\:\: 3 \leq k \leq n-3
\]
\[
\begin{array}{ll}
l=i \neq j & \partial_m \partial_j \partial_k G \\
l = i = j & 2 \partial_m \partial_l \partial_k G
\end{array}
\]
\[
n-2 \leq m \leq n \:\: , \:\:\: 0 \leq k \leq 2
\]
\[
\begin{array}{ll}
l=i \neq j,k & \partial_m \partial_j \partial_k G \\
l = i = j \neq k & 2 \partial_m \partial_l \partial_k G \\
l = i = j = k & 3 \partial_m \partial_l^2 G \: \: .
\end{array}
\]
The entry in the column of $b_{su}$ and the row of $\{i,j,k\}$ is
nonzero only if $s=k$. These possibly nonzero entries are the
following.
\[
n-2 \leq u \leq n \:\: , \:\:\: 3 \leq k \leq n-3
\]
\[
\begin{array}{ll}
s = k & \partial_i \partial_j \partial_u G \: \: .
\end{array}
\]
An easy dimension count shows that we need to prove that there are at most
$6$ relations between the rows of the matrix. Suppose that there are $t$
relations with coefficients
\[
\{ \{ \lambda_{ijk}^r \}\begin{Sb}
0 \leq i,j \leq 2 \\ 0 \leq k \leq n-3
\end{Sb} \}_{1 \leq r \leq t}
\]
between the rows of our matrix. Each relation can be written as a
collection

$3 \leq m \leq n-3$, $0 \leq i \leq 2$
\[
\sum\begin{Sb}
3 \leq k \leq n-3 \\ 0 \leq j \leq 2
\end{Sb} \lambda_{ijk}^r \del_m \del_j \del_k G = 0
\]

$n-2 \leq m \leq n$, $0 \leq i \leq 2$
\begin{equation}
\sum\begin{Sb}
0 \leq k \leq n-3 \\ 0 \leq j \leq 2
\end{Sb} \lambda_{ijk}^r \del_m \del_j \del_k G = 0
\label{secondset}
\end{equation}

$n-2 \leq u \leq n$, $3 \leq k \leq n-3$
\[
\sum\begin{Sb}
0 \leq i,j \leq 2
\end{Sb} \lambda_{ijk}^r \del_i \del_j \del_u G = 0
\]
Each expression $\sum_{0 \leq i,j \leq 2} \lambda_{ijk}^r \del_i \del_j$
defines a hyperplane in $H^0(P, \cO_P(2))$ which contains the polynomials
$\del_u G |_P$. Since we have $3$ independant such polynomials, the vector
space of hyperplanes containing them has dimension $3$. Hence, after a
linear change of coordinates, we can suppose that, for $r \in
\{0, ..., t-3 \}$, we have $\lambda^r_{ijk} = 0$ if $0 \leq i,j \leq 2, 3
\leq k \leq n-3$. The relations (\ref{secondset}) now become

$0 \leq r \leq t-3$, $0 \leq i \leq 2$
\[
\sum\begin{Sb}
0 \leq k \leq 2 \\ 0 \leq j \leq 2
\end{Sb} \lambda_{ijk}^r \del_j \del_k G = 0 \: .
\]
If, for a fixed $r \in \{1, ..., t-3 \}$, the three
relations $\sum\begin{Sb}
0 \leq k \leq 2 \\ 0 \leq j \leq 2
\end{Sb} \lambda_{ijk}^r \del_j \del_k G = 0$, for $0 \leq i \leq 2$,
are not independent, then after a linear change of coordinates, we may
suppose that, for instance, $\lambda_{2jk}^r=0$ for all $j,k \in \{0,1,2
\}$. Since the coefficients $\lambda_{ijk}^r$ are symmetric in $i,j,k$, we
obtain

$0 \leq i \leq 1$
\[
\sum\begin{Sb}
0 \leq k \leq 1 \\ 0 \leq j \leq 1
\end{Sb} \lambda_{ijk}^r \del_j \del_k G = 0 \: .
\]
If $l$ is the line in $P$ obtained as the projectivization of $\langle
\del_0, \del_1 \rangle$, then $\langle \del_{n-2} G, \del_{n-1} G,
\del_n G \rangle |_l$ has dimension at least $2$ and there can be at
most one hyperplane in $H^0(l, \cO_l(2))$ containing $\langle
\del_{n-2} G, \del_{n-1} G, \del_n G \rangle |_l$. In other words, up
to multiplication by a scalar, there is at most one nonzero relation
$\sum\begin{Sb} 0 \leq k \leq 1 \\ 0 \leq j \leq 1
\end{Sb} \lambda_{ijk}^r \del_j \del_k G = 0$. Hence, we can suppose
$\lambda_{1jk}^r=0$ for all $j,k \in \{ 0,1 \}$. Again, by symmetry, we are
reduced to $\lambda_{000}^r \del_0^2 G = 0$ which implies $\lambda_{000}^r
= 0$ because $X$ is smooth. Hence all the $\lambda_{ijk}^r$ are zero.

Therefore, if the $\lambda_{ijk}^r$ are not all zero, the three
relations

$0 \leq i \leq 2$
\[
\sum\begin{Sb}
0 \leq k \leq 2 \\ 0 \leq j \leq 2
\end{Sb} \lambda_{ijk}^r \del_j \del_k G = 0
\]
are independent. If $t-3 \geq 4$, then, after a linear change of
coordinates, for some $r \in \{1, ..., t-3 \}$, one of the above three
relations will be trivial and we are reduced to the previous
case. Therefore $t-3 \leq 3$ and $t \leq 6$.  \hfill $\qed$

\vskip20pt
\begin{proposition}
Suppose that $n \geq 6$. Then $\cP$ has pure dimension equal to the
expected dimension $3n-16$. If $r_P \geq 3$, then $\cP$ is smooth at $P$.
\label{cPsmooth}
\end{proposition}
{\em Proof :} Since the dimension of $\cP_2$ is at most $Min(n-4,5)$ by
Lemma \ref{rP2} and the dimension of every irreducible component of $\cP$
is at least $3n-16$, it is enough to show that for every $P$ such that $r_P
\geq 3$, the space $H^0(P, N_P)$ of infinitesimal deformations of $P$ in
$X$ has dimension $3n-16$.
\vskip10pt

Suppose that $r_P = 3$. As in the proof of Proposition
\ref{rP2inf}, we have an exact sequence
\[0 \lra \cO_P(1)^{\oplus (n-6)} \lra N_{P/X} \lra V_3 \lra 0
\]
where $V_3$ is a locally free sheaf of rank $3$. Since $h^0(P, N_{P/X})
\geq 3n-16$, we have $h^0(P,V_3) \geq 2$. We need to show that
$h^0(P,V_3) = 2$. As in the proof of Proposition \ref{rP2inf} we have
$h^0(P, V_3(-1))=0$ so that, for any line $l \subset P$,
\[
H^0(P, V_3) \inj H^0(l, V_3|_l) \:\: .
\]
Suppose that $h^0(P, V_3) \geq 3$ and let $u_1, u_2, u_3$ be three linearly
independent elements of $H^0(P, V_3)$. By Lemma \ref{linesinP}, the plane
$P$ contains at least one line $l_0$ of second type. It is easily seen that
$V_3|_{l_0} \cong \cO_{l_0}(-1) \oplus \cO_{l_0}(1)^{\oplus 2}$. Therefore
$\langle u_1, u_2, u_3 \rangle |_{l_0}$ generates a subsheaf of the
$\cO_{l_0}(1)^{\oplus 2}$ summand of $V_3 |_{l_0}$ isomorphic to
$\cO_{l_0} \oplus \cO_{l_0}(1)$. The
quotient of $\cO_{l_0}(1)^{\oplus 2}$ by $\cO_{l_0} \oplus \cO_{l_0}(1)$ is
a skyscraper sheaf supported on a point $p$ of $l_0$ (with fiber at $p$
isomorphic to $\bC$). So the images of $u_1, u_2$ and $u_3$ by the
evaluation map at $p$ generate a one-dimensional vector subspace of the
fiber of $V_3$ at $p$. By Lemma \ref{linesinP}, there is a line $l$ of
first type in $P$ which contains $p$. It is easily seen that $V_3|_l \cong
\cO_l^{\oplus 2} \oplus \cO_l(1)$. Restricting $u_1, u_2, u_3$ to $l$ we
see that their images by the evaluation map at $p$ generate a vector
subspace of dimension $\geq 2$ of the fiber of $V_3$ at $p$. Contradiction.

\vskip10pt
Suppose now that $r_P =4$. Then $n \geq 7$ and we have the exact
sequence
\[
0 \lra \cO_P(1)^{\oplus (n-7)} \lra N_{P/X} \lra V_4 \lra 0
\]
where $V_4$ is a locally free sheaf of rank $4$. Since $h^0(P,
N_{P/X}) \geq 3n-16$, we have $h^0(P, V_4) \geq 5$. We need to show
that $h^0(P, V_4) = 5$. As before, $h^0(P, V_4(-1)) = 0$, hence, for
any line $l \subset P$, we have $H^0(P, V_4)
\inj H^0(l, V_4|_l)$. It is easily seen that when $l$ is of first type
$V_4|_l \cong \cO_l^{\oplus 2} \oplus \cO_l(1)^{\oplus 2}$ and when
$l$ is of second type $V_4|_l \cong \cO_l(-1) \oplus
\cO_l(1)^{\oplus 3}$. Thus $h^0(P, V_4) \leq 6$. Suppose that $h^0(P,
V_4) = 6$. Then $H^0(P, V_4) \stackrel{\cong}{\lra} H^0(l, V_4 |_l)$
for every line $l \subset P$.

Suppose first that $P$ contains a line $l_0$ of second type and let $l$
be a line of first type in $P$. We see that $V_4$ is not generated
by its global sections anywhere on $l_0$ whereas $V_4 |_l$ is generated
by its global sections. This gives a contradiction at the point of
intersection of $l$ and $l_0$.

So every line $l$ in $P$ is of first type, $V_4|_l \cong \cO_l^{\oplus 2}
\oplus \cO_l(1)^{\oplus 2}$ and $V_4$ is generated by its global
sections. Let $s$ be a general global section of $V_4$. We claim that $s$
does not vanish at any point of $P$. Indeed, since $V_4$ is generated by its
global sections, for every point $p$ of $P$, the vector space of global
sections of $V_4$ vanishing at $p$ has dimension $2$. Hence the set of all
global sections of $V_4$ vanishing at some point of $P$ has dimension $\leq
2+2 = 4 < 6$. So we have the exact sequence
\[
0 \lra \cO_P \stackrel{s}{\lra} V_4 \lra V \lra 0
\]
where $V$ is a locally free sheaf of rank $3$. Since $V_4$ is generated by
its global sections, so is $V$ and we have $h^0(P, V) = 5$. As before a
general global section $s'$ of $V$ does not vanish anywhere on $P$ and we
have the exact sequence
\[
0 \lra \cO_P \stackrel{s'}{\lra} V \lra V' \lra 0
\] 
where $V'$ is a locally free sheaf of rank $2$. We have $h^0(P, V') = 4$
and $h^0(V'(-1))=h^0(V(-1))=h^0(V_4(-1))=0$. Hence for every line $l
\subset P$, $H^0(P,V') \inj H^0(l,V'|_l)$. Since $V'|_l \cong
\cO_l(1)^{\oplus 2}$, for a nonzero section $s$ of $V'$ the scheme
$Z(s|_l)=Z(s) \cap l$ has length $\leq 1$. The scheme $Z(s)$ is not a line
because $H^0(P,V') \rightarrow H^0(Z(s),V'|_{Z(s)})$ is injective. Hence
for a general line $l \subset P$, $Z(s) \cap l$ is empty. Therefore $Z(s)$
is finite. We compute $c(V') = c(V) = c(V_4) = 1 + 2 \zeta + 4
\zeta^2$. Therefore $Z(s)$ has length $4$. Hence there is a line $l$ such
that $Z(s_l)$ has length $\geq 2$ and this contradicts $length(Z(s_l)) \leq
1$.

\vskip10pt

If $r_P = 5$, consider again the exact sequence of normal sheaves
\[
0 \lra N_{P/X} \lra N_{P/\bP^n} \lra N_{X/ \bP^n} |_P \lra 0
\]
which after tensoring by $\cO_P(-1)$ becomes
\[
0 \lra N_{P/X}(-1) \lra \cO_P^{\oplus (n-2)} \lra \cO_P(2) \lra 0 \: .
\]
Then the map on glabal sections
\[
H^0(P, \cO_P^{\oplus (n-2)}) \lra H^0(P, \cO_P(2))
\]
is onto (see the proof of Proposition \ref{rP2inf}). A fortiori, the map
\[
H^0(P, N_{P/\bP^n})=H^0(P, \cO_P(1)^{\oplus (n-2)}) =
\]
\[
= H^0(P, \cO_P^{\oplus (n-2)}) \otimes H^0(P, \cO_P(1)) \lra H^0(P,
\cO_P(3))= H^0(P, N_{X/ \bP^n} |_P)
\]
is onto and $H^0(P, N_{P/X})$ has dimension $3n-16$.
\hfill
$\qed$

\begin{corollary} If $n \geq 8$, then $\cP$ is irreducible.
\end{corollary}
{\em Proof :} As before, let $G$ be an equation for $X$.  Choose a linear
embedding $\bP^n \inj \bP^{n+1}$. Choose coordinates $\{ x_0, ..., x_n \}$
on $\bP^n$ and coordinates $\{ x_0, ..., x_n, x_{n+1} \}$ on
$\bP^{n+1}$. Let $Y \subset \bP^{n+1}$ be the cubic of equation $G+ x_{n+1}
Q$ where $Q$ is the equation of a general quadric in $\bP^{n+1}$ and let
$\cP_Y \supset \cP$ be the variety of planes in $Y$. Then, by Proposition
\ref{cPsmooth}, the codimension of $\cP$ in $\cP_Y$ is $3$. The singular
locus of $\cP$ is $\cP_2$ (\ref{rP2inf} and \ref{cPsmooth}) which has
codimension at least $4$ in $\cP$ by \ref{rP2} and
\ref{cPsmooth}. Therefore, since $\cP$ is connected (\cite{borcea1} Theorem
4.1 page 33 or \cite{debarremanivel} Th\'eor\`eme 2.1 pages 2-3), it is
sufficient to show that $\cP_Y$ is smooth at a general point of
$\cP_2$. Since $Q$ does not contain a general plane $P \in \cP_2$, the rank
of the dual morphism of $Y$ on $P$ is at least $3$. Hence $\cP_Y$ is smooth
at a general point of $\cP_2$ (\ref{cPsmooth}). \hfill
$\qed$

\begin{lemma}
The dimension of $\cP_3$ is at most $n-2$.
\label{rP3}
\end{lemma}
{\em Proof :} It is enough to show that at any $P$ with
$r_P \leq 3$ the dimension of the tangent space to $\cP_3$ is at most
$n-2$. By \ref{rP2} it is enough to prove this for $r_P = 3$. The
proof of this is very similar to (and simpler than) the proof of Lemma
\ref{rP2}.  \hfill $\qed$

\begin{proposition}
If $n \geq 7$, then $\cP_4$ has pure dimension $2n-9$.
\label{rP4}
\end{proposition}
{\em Proof :} For $n=7$ there is nothing to prove since $\cP$ has
pure dimension $5=3.7-16=2.7-9$ and $\cP = \cP_4$.

Suppose $n \geq 8$. By an easy dimension count, the dimension of every
irreducible component of $\cP_4$ is at least $2n-9$. Since the dimension of
$\cP_3$ is at most $n-2 < 2n-9$ (see \ref{rP3}), for a general element $P$
of any irreducible component of $\cP_4$ we have $r_P = 4$. We first show

\begin{lemma}
Suppose $n \geq 8$. Then the subscheme $\cP_4'$ of $\cP_4$ parametrizing
planes which contain a line of second type has pure dimension $2n-10$.
\label{rP42}
\end{lemma}
{\em Proof :} Again by a dimension count, the dimension of every
irreducible component of $\cP'_4$ is at least $2n - 10$. Let $P$ be an
element of $\cP_4'$. By \ref{rP3}, the scheme $\cP_3 \subset \cP_4'$ has
dimension $\leq n-2 \leq 2n-10$, so we may suppose that $r_P =4$. Let $l$
be the unique line of second type contained in $P$ (see Lemma
\ref{linesinP}). Since the family of lines of second type in $X$ has
dimension $n-3$ (see \cite{CG} Corollary 7.6), it is enough to show that
the space of infinitesimal deformations of $P$ in $X$ which contain $l$ has
dimension $n-7$.

Consider the exact sequence of sheaves
\[ 0 \lra N_{P/X}(-1) \lra N_{P/X} \lra N_{P/X} |_l \lra 0 \]
with associated cohomology sequence
\[
0 \lra H^0(P, N_{P/X}(-1)) \lra H^0(P, N_{P/X}) \lra H^0(P, N_{P/X} |_l)
\lra H^1(P, N_{P/X}(-1)) \lra ...
\]
The space of infinitesimal deformations of $P$ in $X$ which contain $l$ can
be identified with the kernel of the homomorphism $H^0(P, N_{P/X}) \lra
H^0(P, N_{P/X} |_l) $ which, by the above sequence, can be
identified with $H^0(P, N_{P/X}(-1))$. Recall the exact sequence
\[
0 \lra N_{P/X}(-1) \lra \cO_P^{\oplus (n-2)} \lra \cO_P(2) \lra 0
\]
where the map $\cO_P^{\oplus (n-2)} \lra \cO_P(2)$ is given by
multiplication by $\del_3 G, ..., \del_n G$ (see the proof of
\ref{rP2inf}). It immediately follows that $h^0(P, N_{P/X}(-1)) = n-7$ if
and only if $r_P=4$.  \hfill $\qed$ \vskip20pt

Note that containing a line of second type imposes at most one
condition on planes $P$ with $r_P \leq 4$. Therefore Proposition
\ref{rP4} follows from Lemma \ref{rP42}.
\hfill
$\qed$

\section{Resolving the indeterminacies of the rational involution on
$S_l$}
\label{secSl'}

A good generalization of the Prym construction for cubic threefolds to
cubic hypersurfaces of higher dimension would be to realize the cohomology
of $X$ as the anti-invariant part of the cohomology of $S_l$ for the
involution exchanging two lines whenever they are in the same fiber of
$\pi$. However, this is only a rational involution and we need to resolve
its indeterminacies. This involution is not well-defined exactly at the
lines $l'$ such that $\pi^{-1}(\pi(l')) \subset X_l$, i.e., the plane $L'
\subset \bP^n$ corresponding to $\pi(l')$ is contained in $X$. Let $T_l
\subset Q_l \subset \bP^{n-2}$ be the variety parametrizing the planes in
${\Bbb P}^n$ which contain $l$ and are contained in $X$ (equivalently, the
variety $T_l$ parametrizes the fibers of $\pi$ which are contained in
$X_l$). Recall that $X_l \subset {\Bbb P}^n_l$ is the divisor of zeros of
$s \in H^0(\bP E, \cO_{\bP E}(2) \otimes \pi^* \cO_{\bP^{n-2}}(1)) =
H^0(\bP^{n-2}, \pi_*(\cO_{\bP E}(2)) \otimes \cO_{\bP^{n-2}}(1)) =
H^0(\bP^{n-2}, Sym^2E^* \otimes \cO_{\bP^{n-2}}(1))$. Since $E \cong
\cO_{{\Bbb P}^{n-2}}(-1) \oplus \cO_{\bP^{n-2}}^{\oplus 2}$, we have
$Sym^2E^* \otimes \cO_{\bP^{n-2}}(1) \cong \cO_{\bP^{n-2}}(3) \oplus
\cO_{\bP^{n-2}}(2)^{\oplus 2} \oplus \cO_{\bP^{n-2}}(1)^{\oplus 3}$. The
variety $T_l$ is the locus of common zeros of all the components of $s$ in
the above direct sum decomposition. Therefore $T_l$ is the scheme-theoretic
intersection of three hyperplanes, two quadrics and one cubic in
$\bP^{n-2}$. We have

\begin{lemma}
There is a Zariski-dense open subset of $F$ parametrizing lines $l$ such
that $l$ is of first type and $r_P=5$ for every plane $P$ in $X$ containing
$l$. For $l$ in this Zariski-dense open subset, the variety $T_l$ is the
smooth complete intersection of the six hypersurfaces obtained as the zero
loci of the components of $s$ in the direct sum decomposition $Sym^2E^*
\otimes \cO_{\bP^{n-2}}(1) \cong \cO_{\bP^{n-2}}(3) \oplus
\cO_{\bP^{n-2}}(2)^{\oplus 2} \oplus \cO_{\bP^{n-2}}(1)^{\oplus 3}$.
\label{Tlsmooth}
\end{lemma}
{\em Proof :} The first part of the lemma follows from Proposition
\ref{rP4}. For the second part we need to show that $T_l$ is smooth of the
expected dimension $n-8$. In other words, for any plane $P$ containing $l$,
the space of infinitesimal deformations of $P$ in $X$ containing $l$ has
dimension $n-8$. The proof of this is similar to the proof of
Lemma \ref{rP42}.  \hfill \qed \vskip20pt

\begin{definition}
Let $U_0$ be the subvariety of $F$ parametrizing lines $l$ such that $l$ is
of first type, is not contained in $\bT_{l'}$ for any line $l'$ of second
type and every plane containing $l$ is an element of $\cP \setminus \cP_4$.
\end{definition}

By Lemmas \ref{lemTl'} and \ref{Tlsmooth}, the variety $U_0$ is an open
dense subvariety of $F$. Suppose $l \in U_0$. By Lemmas
\ref{Dlsmooth1}, \ref{lemSlsm} and \ref{Tlsmooth}, the varieties $S_l$
and $T_l$ are smooth of the expected dimensions $n-3$ and $n-8$
respectively. Let $X'_l \subset {\bP^n_l}'$ be the blow ups of $X_l
\subset \bP^n_l$ along $\pi^{-1}(T_l)$ and let ${\bP^{n-2}}'$ be the
blow up of $\bP^{n-2}$ along $T_l$. Then we have morphisms
\[ \begin{array}{crl}
X'_l & \subset & {\bP^n_l}' \\
 & \pi_X' \searrow & \downarrow \pi' \\
 & & {\bP^{n-2}}'
\end{array} \]
where $\pi' : {\bP^n_l}' \lra {\bP^{n-2}}'$ is again a
$\bP^2$-bundle. Since $T_l$ is the zero locus of $s \in H^0(\bP^{n-2},
\pi_* \cO_{\bP E}(2) \otimes \cO_{\bP^{n-2}}(1))$, we have
$N_{T_l/\bP^{n-2}} \cong \pi_* \cO_{\bP E}(2) \otimes
\cO_{\bP^{n-2}}(1) |_{T_l}$. Therefore, the exceptional divisor $E'$
of ${\bP^{n-2}}' \lra \bP^{n-2}$ is a $\bP^5$-bundle over $T_l$ whose
fiber at a point $t \in T_l$ corresponding to the plane $P_t \subset
X_l$ is $| \cO_{P_t}(2) |$. We have

\begin{lemma}
Suppose that $l \in U_0$. For all $t \in T_l$, the restriction of $\pi_X' :
X_l' \lra {\bP^{n-2}}'$ to $| \cO_{P_t}(2) | \subset {\bP^{n-2}}'$ is the
universal conic on $| \cO_{P_t}(2) |$. In particular, the fibers of $\pi_X'
: X'_l \lra {\bP^{n-2}}'$ are always one-dimensional.
\label{fibonedim}
\end{lemma}

{\em Proof :} The restriction of $\pi'$ to the inverse image of a point $t
\in T_l$ is the second projection $P_t \times | \cO_{P_t}(2) | \lra |
\cO_{P_t}(2) |$. Let $N_{X,p}$ be the normal space in $X_l$ to
$\pi^{-1}(T_l)$ at $p \in P_t$ and let $\rho_t : P_t \lra | \cO_{P_t}(2)|^*
\cong \bP^5$ be the map which to $p \in P_t$ associates $\bP N_{X,p} \in |
\cO_{P_t}(2) |^*$. For $n \in | \cO_{P_t}(2) |$, the fiber of $\pi'_X$ at
$(t,n) \in E'$ is equal to $\rho_t^{-1}(\rho_t(P_t) \cap H_n)$ where $H_n$
is the hyperplane in $| \cO_{P_t}(2) |^*$ corresponding to $n$. It is
immediately seen that $\rho_t$ is induced by the dual morphism $\delta$ of
$X$. Hence, since $r_{P_t} = 5$, the map $\rho_t$ is the Veronese morphism
$P_t \lra | \cO_{P_t}(2) |^*$. Hence $\rho_t^{-1}(\rho_t(P_t) \cap H_n)$ is
the conic in $P_t$ corresponding to $n$. \hfill $\qed$

\vskip20pt
It follows from lemma \ref{fibonedim} that if we let $S'_l$ be the
variety parametrizing lines in the fibers of $\pi_X' : X'_l \lra
{\bP^{n-2}}'$, then there is a well-defined involution $i_l : S'_l \lra
S_l'$ which sends $l'$ to $l''$ when $l'+l''$ is a fiber of $X'_l \lra
{\bP^{n-2}}'$. Sending a line in a fiber of $\pi'_X$ to its image in $X_l$
defines a morphism $S_l' \lra S_l$. Let $\cP_l \rightarrow T_l$ be the
family of planes in $X$ containing $l$, then the inverse image of $T_l$ in
$S_l$ by the morphism $S_l \rightarrow Q_l$ is the projective bundle
$\cP_l^*$ of lines in the fibers of $\cP_l \rightarrow T_l$. We have

\begin{proposition}
Suppose that $l \in U_0$. The morphism $S_l' \lra S_l$ is the blow up of
$S_l$ along $\cP^*_l$. In particular, the variety $S_l'$ is smooth. The
fixed point locus $R_l'$ of $i_l$ in $S_l'$ is a smooth subvariety of
codimension $2$ of $S_l'$. The projective bundle $\bP(N_{R_l'/S_l'}) \lra
R_l'$ is isomorphic to the family of lines in the fibers of $\pi_X'$
parametrized by $R_l'$.
\label{everythingsmooth}
\end{proposition}
{\em Proof :} In Lemma \ref{lemSlsm}, we saw that $S_l$ can be
identified with the closure of the subvariety $G(2,n+1) \times G(3,n+1)$
parametrizing pairs $(l',L')$ of a line and a plane such that $l \neq l'$
and $l \cup l' \subset L'$. In the same way, we see that $S_l'$ can be
identified with the closure of the subvariety of $G(2,n+1) \times G(2,n+1)
\times G(3,n+1)$ parametrizing triples $(l',l'',L')$ such that $L' \cap X
\supset l \cup l' \cup l''$ and $l, l', l''$ are distinct. Furthermore, the
morphism $S'_l \rightarrow S_l$ is the restriction of the projection to the
second and third factors of $G(2,n+1) \times G(2,n+1) \times
G(3,n+1)$. Again as in the proof of Lemma \ref{lemSlsm} we see that $S_l'$
is smooth. Blowing up $\cP_l^*$ and its inverse image in $S_l'$ we obtain
the commutative diagram
\[
\begin{array}{ccc}
\tS_l' & \lra & \tS_l \\
\downarrow & & \downarrow \\
S_l' & \lra & S_l .
\end{array}
\]
Since the inverse image of $\cP_l^*$ is a divisor in $S_l'$, the blow up
morphism $\tS_l' \rightarrow S_l'$ is an isomorphism. The morphism $S_l'
\rightarrow \tS_l$ thus obtained is a birational morphism of smooth
varieties with constant fiber dimension hence it is an isomorphism. This
proves the first part of the Proposition.

Now let $\Delta$ be the diagonal of $G(2,n+1) \times G(2,n+1)$. Then the
variety $R_l'$ is identified with $S'_l \cap \left( \Delta \times G(3,n+1)
\right)$. One now computes the tangent space to $R_l'$ as in the proof of
Lemma \ref{lemSlsm} and see that $N_{R_l'/S_l'}$ is isomorphic to $I^*
\otimes J/I$ where $I$ is the restriction of the universal bundle on
$G(2,n+1)$ and $J$ is the restriction of the universal bundle on
$G(3,n+1)$. Therefore $\bP(N_{R_l'/S_l'})$ is isomorphic to $\bP(I)$ which
is the family of lines in the fibers of $\pi_X'$ parametrized by $R_l'$.
\hfill \qed \vskip15pt

Let $Q_l'$ be the blow up of $Q_l$ along $T_l$. Sending a line $l \in S_l'$
to the fiber of $X'_l \lra {\bP^{n-2}}'$ which contains it defines a finite
morphism $S'_l \lra Q_l'$ of degree $2$ with ramification locus
$R_l'$. Blowing up $R_l'$ in $Q_l'$ and $S_l'$ we obtain the morphism
$S_l'' \lra Q_l''$. We have
\begin{proposition} \label{propQ'''}
The variety $R_l'$ is an ordinary double locus for $Q_l'$. In particular,
$Q_l''$ is smooth and (by \ref{everythingsmooth}) the projectivization
$\bP(C_{R_l'/Q_l''})$ of the normal cone to $R_l'$ in $Q_l'$ is isomorphic
to $\bP(N_{R_l'/S_l'})$.
\end{proposition}

{\em Proof :} The fact that $R_l \setminus T_l$ is an ordinary double locus
for $Q_l \setminus T_l$ can be proved, for instance, by intersecting $Q_l$
with a general plane through a point $p$ of $R_l \setminus T_l$. The
resulting curve has an ordinary double point at $p$ by \cite{B2}
Proposition 1.2 page 321. At a point $q$ of the exceptional divisor of
$R_l' \rightarrow R_l$, locally trivialize the pull-back of
$E=\cO_{\bP^{n-2}}^{\oplus 2} \oplus \cO_{\bP^{n-2}}(-1)$ to obtain a
morphism from a neighborhood $U$ of $q$ to $|\cO_{\bP^2}(2)|$. It easily
follows from \ref{Tlsmooth} and \ref{fibonedim} that this morphism is
dominant and the restriction of $X_l' \rightarrow {\bP^{n-2}}'$ to $U$ is
the inverse image of the universal conic on $|\cO_{\bP^2}(2)|$. The
assertion of the Proposition now follows from the corresponding fact for
the cubic fourfold parametrizing singular conics in $\cP^2$. \hfill \qed
\vskip15pt

\section{The main Theorem}
\label{secmaintheorem}

Let $L_l \lra S'_l$ and $\oL_l \lra S_l$ be the families of lines in the
fibers of $\pi'_X$ and $\pi_X$ respectively. The blow-up morphism
$\epsilon_2 : X_l' \lra X_l$ defines a morphism $L_l \lra \oL_l$ which fits
into the commutative diagram
\[ \begin{array}{rclcc}
X_l' & \stackrel{\epsilon_2}{\lra} & X_l & \stackrel{\epsilon_1}{\lra}
& X \\
\rho \uparrow & & \uparrow \orho & & \\
L_l & \lra & \oL_l & & \\
p \downarrow & & \downarrow \op & & \\
S_l' & \lra & S_l & &  
\end{array}
\]
where the squares are Cartesian. Put $q = \epsilon_1 \epsilon_2 \rho$ and
let $\psi' = q_* p^* : H^{n-3}(S'_l, \bZ) \lra H^{n-1}(X, \bZ)$ and $\psi =
(\epsilon_1 \orho)_* \op^*: H^{n-3}(S_l, \bZ) \lra H^{n-1}(X, \bZ)$ be the
Abel-Jacobi maps. The map $\psi$ is the composition of $\psi'$ with the
inclusion $H^{n-3}(S_l, \bZ) \inj H^{n-3}(S'_l, \bZ)$ because the bottom
(or top) square above is Cartesian. We have

\begin{theorem}
The maps $\psi : H^{n-3}(S_l, \bZ) \lra H^{n-1}(X, \bZ)$ and $\psi' : 
H^{n-3}(S'_l, \bZ) \lra H^{n-1}(X, \bZ)$ are onto.
\label{thmpsionto}
\end{theorem}

{\em Proof :} Consider the rational map $Q_l' \lra X_l'$ which to the
singular conic $l'+l''$ associates the point of intersection $l' \cap
l''$. An easy local computation shows that the closure of the image of this
map is smooth, hence, by a reasoning analogous to the proof of Proposition
\ref{everythingsmooth}, it can be identified with $Q_l''$. Let $\epsilon_3
: X_l'' \lra X_l'$ be the blow up of $X_l'$ along $Q_l''$ and, for each $i$
($1 \leq i \leq 3$), let $E_i$ be the exceptional divisor of the blow up
map $\epsilon_i$. Then we have a factorization
\[ \begin{array}{crl}
 & & X_l'' \\
 & \tq \nearrow & \downarrow \epsilon_3\\
L_l & \stackrel{\rho}{\lra} & X_l'
\end{array} \]
so that $\psi' = q_*p^* = {\epsilon_1}_* {\epsilon_2}_* \rho_* p^* =
{\epsilon_1}_* {\epsilon_2}_* {\epsilon_3}_* \tq_* p^* $. Note that
$\tq$ is an embedding so that we can, and will, identify $L_l$ with
$\tq(L_l)$. Put $U_l = X''_l \setminus (E_3 \cup L_l) = X'_l
\setminus \rho(L_l)$.  Let $m_l : U_l \lra X''_l$ be the inclusion. We
have the spectral sequence
\[
E^{p,q}_2 = H^p(X''_l, R^q{m_l}_* \bZ_{U_l}) \implies H^{p+q}(U_l, \bZ)
\]
and by \cite{deligne}, $\S$3.1, we have $R^0{m_l}_* \bZ_{U_l} = \bZ_{X''_l},
R^1{m_l}_* \bZ_{U_l} = \bZ_{E_3} \oplus \bZ_{L_l}, R^2{m_l}_* \bZ_{U_l} =
\bZ_{E_3 \cap L_l}$ and $R^q {m_l}_* \bZ_{U_l} = 0$ for $q > 2$. Note that
$E_3 \cap L_l \cong S_l''$.

Therefore
\[ \begin{array}{c}
E^{p,0}_2 = H^p(X''_l, \bZ) \: , \\
E^{p,1}_2 = H^p(X''_l, \bZ_{E_3} \oplus \bZ_{L_l}) = H^p(L_l, \bZ) \oplus
H^p(E_3, \bZ) \: , \\
E^{p,2}_2 = H^p(X''_l, \bZ_{S_l''} ) = H^p(S''_l, \bZ) \: , \\
E^{p,q}_2 = 0 \hbox{ for } q > 2
\end{array}
\]
So the $E_2^{.,.}$ complex is
\[ 0 \lra H^{p-2}(S''_l, \bZ) \lra H^p(L_l, \bZ) \oplus H^p(E_3, \bZ)
\lra H^{p+2}(X''_l, \bZ) \lra 0
\]
where the maps are obtained by
Poincar\'e Duality from the natural push-forwards on homology induced
by the inclusions. We have (see, for instance, \cite{B2}, 0.1.3, page 312)

\begin{eqnarray}H^{p+2}(X''_l, \bZ) \cong H^{p+2}(X'_l, \bZ) \oplus
H^p(Q''_l, \bZ) \:,
\label{cohomXlll} \\
H^{p+2}(X'_l, \bZ) \cong
H^{p+2}(X_l, \bZ) \oplus \left( \oplus\begin{Sb}
p-6 \leq i \leq p \\ i \equiv p [2]
\end{Sb} H^i(\pi^{-1}(T_l), \bZ) \right) \: ,
\label{cohomXll} \\
H^{p+2}(X_l, \bZ) \cong H^{p+2}(X, \bZ) \oplus \left( \oplus\begin{Sb}
p-2(n-4) \leq i \leq p \\ i \equiv p [2]
\end{Sb} H^i(l, \bZ) \right) \label{cohomXl}
\end{eqnarray}
and
\begin{eqnarray}
H^{p-2}(S''_l, \bZ) \cong H^{p-2}(S'_l, \bZ) \oplus H^{p-4}(R'_l, \bZ).
\end{eqnarray}
Since $E_3$ and $L_l$ are $\bP^1$-bundles over $Q''_l$ and $S'_l$
respectively,
\begin{eqnarray} H^p(E_3, \bZ) \cong H^p(Q''_l, \bZ) \oplus
H^{p-2}(Q''_l, \bZ) \end{eqnarray}
and
\begin{eqnarray} H^p(L_l, \bZ) \cong H^p(S'_l, \bZ) \oplus
H^{p-2}(S'_l, \bZ) \: . \label{cohomLl} \end{eqnarray}

The map $\psi'$ is the composition of the
inclusion $H^{n-3}(S'_l, \bZ) \inj H^{n-3}(L_l, \bZ)$ obtained from
(\ref{cohomLl}) with the differential $E^{n-3,1}_2 \lra E^{n-1,0}_2$
and the projection $H^{n-1}(X_l'', \bZ) \surj H^{n-1}(X, \bZ)$
obtained from (\ref{cohomXlll}), (\ref{cohomXll}) and
(\ref{cohomXl}). We first study the cokernel of the differential
$E^{n-3,1}_2 \lra E^{n-1,0}_2$.

By \cite{deligne}, 3.2.13, the differentials $E^{p,q}_3 \lra E^{p+3,
q-2}_3$ are zero. Therefore $E^{.,.}_{\infty} = E^{.,.}_3$ and, in
particular,
\[ Coker \left( H^{n-3}(L_l, \bZ) \oplus H^{n-3}(E_3, \bZ)
\lra H^{n-1}(X''_l, \bZ) \right) = \]
\[ = Coker \left( E_2^{n-3,1} \lra E_2^{n-1,0} \right) \]
\[ = E_3^{n-1,0} = E_{\infty}^{n-1,0} = Gr^{n-1}(H^{n-1}(U_l, \bZ))
\: . \]
This is the image of $H^{n-1}(X''_l, \bZ)$ in $H^{n-1}(U_l, \bZ)$ and,
by \cite{deligne} 3.2.17, it is the piece $W_{n-1}(H^{n-1}(U_l, \bZ))$
of weight $n-1$ of the mixed Hodge structure on $H^{n-1}(U_l, \bZ)$.

Define $V_l := {\bP^{n-2}}' \setminus Q_l'$. The fibers of the
conic-bundle $U_l \lra V_l$ are all smooth, hence
\[H^{n-1}(U_l, \bZ) \cong H^{n-3}(V_l, \bZ) \oplus H^{n-1}(V_l, \bZ)
\]
\begin{claim} Under this isomorphism, the space
$W_{n-1}(H^{n-1}(U_l, \bZ))$ is isomorphic to \linebreak
$W_{n-3}(H^{n-3}(V_l, \bZ)) \oplus W_{n-1}(H^{n-1}(V_l, \bZ))$.
\end{claim}
To prove this it is sufficient to show that the maps $H^{n-1}(V_l, \bZ)
\lra H^{n-1}(U_l, \bZ)$ and $H^{n-3}(V_l, \bZ) \lra H^{n-1}(U_l, \bZ)$ are
morphisms of mixed Hodge structures of type $(0,0)$ and $(1,1)$
respectively.

By \cite{deligne} pages 37-38, the pull-backs on cohomology
$H^{n-3}(V_l, \bZ) \lra H^{n-3}(U_l, \bZ)$ and $H^{n-1}(V_l, \bZ) \lra
H^{n-1}(U_l, \bZ)$ are morphisms of mixed Hodge structures of type
$(0,0)$. To see that the map $H^{n-3}(V_l, \bZ) \lra H^{n-1}(U_l,
\bZ)$ is a morphism of mixed Hodge structures of type $(1,1)$
choose a bisection $B$ of the conic bundle $U_l
\lra V_l$ and let $\eta$ be a half of the cohomology class of
$B$. Then the map
\[ H^{n-3}(V_l, \bZ) \lra H^{n-1}(U_l, \bZ) \]
is the composition of pull-back
\[ H^{n-3}(V_l, \bZ) \lra H^{n-3}(U_l, \bZ) \]
with cup-product with $\eta$
\[ H^{n-3}(U_l, \bZ) \lra H^{n-1}(U_l, \bZ) \: .\]
The class $2 \eta$ is the restriction to $U_l$ of the cohomology class
of the closure of $B$ in $X'_l$. Therefore $2 \eta$ is in the image
of
\[ H^2(X'_l, \bZ) \lra H^2(U_l, \bZ) \]
and hence has pure weight $2$ and Hodge type $(1,1)$. Therefore $\eta$ has
pure weight $2$ and Hodge type $(1,1)$ in the mixed Hodge structure on
$H^2(U_l, \bZ)$, the map $H^{n-3}(V_l, \bZ) \lra H^{n-1}(U_l, \bZ)$ is a
morphism of mixed Hodge structures of type $(1,1)$ and sends
$W_{n-3}(H^{n-3}(V_l, \bZ))$ into $W_{n-1}( H^{n-1}(U_l, \bZ))$.
\vskip20pt
 
We now determine $W_{n-3}(H^{n-3}(V_l, \bZ)) \oplus
W_{n-1}(H^{n-1}(V_l, \bZ))$. In the following we let $p$ be equal to
$n-3$ or $n-1$.

Let ${\bP^{n-2}}'' \lra {\bP^{n-2}}'$ be the blow up of ${\bP^{n-2}}'$
along $R_l'$ with exceptional divisor $E''$ and identify $Q_l''$ with
its image in ${\bP^{n-2}}''$. Then $V_l = {\bP^{n-2}}'' \setminus
\left( E'' \cup Q_l'' \right)$ and the divisors $E''$ and $Q_l''$ are
smooth and meet transversally. Therefore $W_p(H^p(V_l, \bZ))$ is the
image of $H^p({\bP^{n-2}}'', \bZ)$ in $H^p(V_l, \bZ)$, i.e., it is
isomorphic to the cokernel of the map

\[ H^{p-2}(Q_l'', \bZ) \oplus
H^{p-2}(E'', \bZ) \lra H^p({\bP^{n-2}}'', \bZ) \]
obtained by Poincar\'e Duality from push-forward on homology. Since
$E''$ is a $\bP^2$-bundle over $R_l'$, we have
\begin{eqnarray}
H^{p-2}(E'', \bZ) \cong H^{p-2}(R_l', \bZ) \oplus H^{p-4}(R_l', \bZ)
\oplus H^{p-6}(R_l', \bZ) \: ,
\end{eqnarray}
By e.g. \cite{B2} 0.1.3, we have the isomorphism
\[
H^p({\bP^{n-2}}'', \bZ) \cong H^p({\bP^{n-2}}', \bZ) \oplus
H^{p-2}(R_l', \bZ) \oplus H^{p-4}(R_l', \bZ) \: .
\]
Under the map $H^{p-2}(E'', \bZ) \lra H^p({\bP^{n-2}}'', \bZ)$ above, the
summand $H^{p-2}(R_l', \bZ) \oplus H^{p-4}(R_l', \bZ)$ in $H^{p-2}(E'',
\bZ)$ maps isomorphically onto the same summand in $H^p({\bP^{n-2}}'',
\bZ)$. Therefore $W_p(H^p(V_l, \bZ))$ is a quotient of
$H^p({\bP^{n-2}}',\bZ)$.

The summand $H^{p-6}(R_l', \bZ)$ in $H^{p-2}(E'', \bZ)$ maps into the
summand $H^p({\bP^{n-2}}', \bZ)$ of $H^p({\bP^{n-2}}'', \bZ)$, the map
$H^{p-6}(R_l', \bZ) \lra H^{p}({\bP^{n-2}}', \bZ)$ being again
obtained by Poincar\'e Duality from push-forward on homology. Since
the degree of $R_l$ in $\bP^{n-2}$ is $16$, the image of the
composition of $H^{p-6}(R_l', \bZ) \inj H^{p}({\bP^{n-2}}', \bZ)$ with
the isomorphism
\[
H^p({\bP^{n-2}}', \bZ) \cong H^p(\bP^{n-2}, \bZ) \oplus
\left( \oplus\begin{Sb}
p-10 \leq i \leq p-2 \\ i \equiv p [2]
\end{Sb} H^i(T_l, \bZ) \right)
\]
contains an element whose component in the summand $H^p(\bP^{n-2},
\bZ)$ is $16$ times a generator of $H^p(\bP^{n-2}, \bZ)$.

Since the degree of $Q_l$ is $5$, the image of the composition of the
direct sum embedding
\[
H^{p-2}(Q_l'', \bZ) \inj H^{p-2}(E'', \bZ) \oplus H^{p-2}(Q_l'', \bZ)
\]
with the map
\[
H^{p-2}(E'', \bZ) \oplus H^{p-2}(Q_l'', \bZ) \lra H^p({\bP^{n-2}}'', \bZ)
\]
contains an element whose component in the summand $H^p({\bP^{n-2}},
\bZ)$ is $5$ times a generator of $H^{p}(\bP^{n-2}, \bZ)$. Since $16$
and $5$ are coprime, we deduce that the image of $H^{p-2}(E'', \bZ)
\oplus H^{p-2}(Q_l'', \bZ)$ in $H^{p}({\bP^{n-2}}'', \bZ)$ contains
an element whose component in the summand $H^{p}({\bP^{n-2}}, \bZ)$ is
a generator of $H^p({\bP^{n-2}}, \bZ)$.
\vskip10pt

{\em So far we have obtained that $W_{p}(H^p(V_l, \bZ))$ is a
quotient of}
\[
\oplus\begin{Sb}
p-10 \leq i \leq p-2 \\ i \equiv p [2]
\end{Sb} H^i(T_l, \bZ) \subset H^p({\bP^{n-2}}'', \bZ)\: .
\]
It is now easily seen that $\left( \oplus\begin{Sb}
n-11 \leq i \leq n-3 \\ i \equiv n-1 [2]
\end{Sb} H^i(T_l, \bZ) \right) \oplus \left( \oplus\begin{Sb}
n-13 \leq i \leq n-5 \\ i \equiv n-1 [2]
\end{Sb} H^i(T_l, \bZ) \right)$
maps into the summand
\[
\oplus\begin{Sb}
n-9 \leq i \leq n-3 \\ i \equiv n-3 [2]
\end{Sb} H^i(\pi^{-1}(T_l), \bZ)
\]
of $H^{n-1}(X_l'', \bZ)$. Therefore $W_{n-1}(H^{n-1}(U_l, \bZ))=
W_{n-3}(H^{n-3}(V_l, \bZ)) \oplus W_{n-1}(H^{n-1}(V_l, \bZ))$ is a
subquotient of
\[
\oplus\begin{Sb}
n-9 \leq i \leq n-3 \\ i \equiv n-3 [2]
\end{Sb} H^i(\pi^{-1}(T_l), \bZ) \subset H^{n-1}(X_l'', \bZ) =
\]
\[
H^{n-1}(X_l, \bZ) \oplus H^{n-3}(Q_l'', \bZ) \oplus \left(
\oplus\begin{Sb} n-9 \leq i \leq n-3 \\ i \equiv n-3 [2] \end{Sb}
H^i(\pi^{-1}(T_l), \bZ) \right)
\]
and the map
\[
H^{n-3}(L_l, \bZ)\oplus H^{n-3}(E_3, \bZ) \lra H^{n-1}(X_l, \bZ)
\]
is onto. So, in particular, we have proved
\begin{claim}
The map
\[
H^{n-3}(L_l, \bZ) \oplus H^{n-3}(E_3, \bZ) \lra H^{n-1}(X, \bZ)
\]
is onto.
\end{claim}

Since $E_3$ is the exceptional divisor of the blow up $X_l'' \lra
X_l'$, the image of
\[
H^{n-3}(E_3, \bZ) \lra H^{n-1}(X, \bZ)
\]
is equal to the image of
\[
H^{n-5}(Q_l'', \bZ) \lra H^{n-1}(X, \bZ) \: .
\]
We will prove that the image of this map is algebraic. Since
$H^{n-1}(X,\bZ)$ is torsion-free, it is enough to
prove this after tensoring with $\bQ$. Since, by Poincar\'e Duality,
$H^{n-5}(Q_l'', \bQ) \cong H^{n-1}(Q_l'', \bQ)^*$, we first determine
$H^{n-1}(Q_l'', \bQ)$. For this we use the spectral sequence
\[
E^{p,q}_2 = H^p({\bP^{n-2}}'', R^q{u}_* \bZ) \implies H^{p+q}(W, \bZ)
\]
where $W := \bP^{n-2} \setminus Q_l = {\bP^{n-2}}'' \setminus \left(
\tE' \cup E'' \cup Q_l'' \right)$ with $\tE'$ the proper transform of
$E'$ in ${\bP^{n-2}}''$ and $u : W \inj {\bP^{n-2}}''$ is the
inclusion. Recall that such a spectral sequence degenerates at $E_3$
(\cite{deligne}, 3.2.13). By \cite{goreskymcpherson} pages 23-24, we
have $H^i(W, \bZ) = 0$ for $i > dim(W) = n-2$. Therefore we obtain the
following exact sequence from the spectral sequence
\begin{eqnarray}
H^{n-5}(\tE' \cap E'' \cap Q_l'', \bZ) \stackrel{d_{n-3}}{\lra}
H^{n-3}(\tE' \cap E'', \bZ) \oplus H^{n-3}(\tE' \cap Q_l'', \bZ)
\oplus H^{n-3}(E'' \cap Q_l'', \bZ) \stackrel{d_{n-1}}{\lra} \nonumber
\end{eqnarray}
\begin{eqnarray}
\stackrel{d_{n-1}}{\lra} H^{n-1}(\tE', \bZ) \oplus H^{n-1}(E'', \bZ)
\oplus H^{n-1}(Q_l'', \bZ) \stackrel{d_{n+1}}{\lra}
H^{n+1}({\bP^{n-2}}'', \bZ) \lra 0 \:\: .
\label{Qlcoh}
\end{eqnarray}
We have

\begin{lemma} The varieties whose cohomologies appear in
sequence (\ref{Qlcoh}) are described as follows.
\label{descrinters}
\begin{description}

\item[$\bold \tE' \cap E'' \cap Q_l''$] $\bP^1$-bundle over $\cV_l$
where $\cV_l := E' \cap R_l'$. The variety $\cV_l$ is a $\bP^2$-bundle
over $T_l$ and each of its fibers over $T_l$ embeds into the
corresponding fiber of $E'$ as the Veronese surface. Hence
\[
H^{n-5}(\tE' \cap E'' \cap Q_l'', \bZ) \cong H^{n-5}(\cV_l, \bZ)
\oplus H^{n-7}(\cV_l, \bZ)
\]
and
\[
H^i(\cV_l, \bZ) \cong H^i(T_l, \bZ) \oplus H^{i-2}(T_l, \bZ) \oplus
H^{i-4}(T_l, \bZ) \:\: .
\]

\item[$\bold T_l'' := \tE' \cap Q_l''$] bundle over $T_l$ with fibers
isomorphic to the blow up $\hS^2 \bP^2$ of the symmetric square $S^2
\bP^2$ of $\bP^2$ along the diagonal of $S^2 \bP^2$. A fiber of $\tE'
\cap E'' \cap Q_l''$ embeds into the corresponding fiber of $\tE' \cap
Q_l''$ as the exceptional divisor of the blow up $\hS^2 \bP^2 \lra S^2
\bP^2$. We have
\[
H^{n-3}(T_l'', \bZ) \cong H^{n-3}(T_l, \bZ) \oplus
H^{n-5}(T_l, \bZ) \oplus H^{n-7}(T_l, \bZ)^{\oplus 2} \oplus
H^{n-9}(T_l, \bZ) \oplus H^{n-11}(T_l, \bZ) \oplus H^{n-5}(\cV_l, \bZ)
\]
\[
\cong H^{n-3}(T_l, \bZ) \oplus H^{n-5}(T_l, \bZ) \oplus H^{n-7}(T_l,
\bZ) \oplus H^{n-7}(\cV_l, \bZ) \oplus H^{n-5}(\cV_l, \bZ)
\]
and, under $d_{n-3}$, the summand $H^{n-7}(\cV_l,
\bZ) \oplus H^{n-5}(\cV_l, \bZ)$ in $H^{n-5}(\tE' \cap E'' \cap Q_l'',
\bZ)$ maps into the same summand in $H^{n-3}(T_l'', \bZ)$.

\item[$\bold E'' \cap Q_l''$] $\bP^1$-bundle over $R_l'$. Hence
\[
H^{n-3}(E'' \cap Q_l'', \bZ) \cong H^{n-3}(R_l', \bZ) \oplus
H^{n-5}(R_l', \bZ) \:\: .
\]

\item[$\bold \tE' \cap E''$] $\bP^2$-bundle over $\cV_l$ which
contains $\tE' \cap E'' \cap Q_l''$ as a conic-bundle over $\cV_l$. We
have
\[
H^{n-3}(\tE' \cap E'', \bZ) \cong H^{n-3}(\cV_l, \bZ) \oplus
H^{n-5}(\cV_l, \bZ) \oplus H^{n-7}(\cV_l, \bZ)\:\: .
\]

\item[$\bold \tE'$] the blow up of $E'$ along $\cV_l$, i.e., bundle
over $T_l$ with fibers isomorphic to the blow up of $\bP^5$ along the
Veronese surface. This contains $\tE' \cap E''$ as its exceptional
divisor. Hence
\[
H^{n-1}(\tE', \bZ) \cong H^{n-3}(\cV_l, \bZ) \oplus H^{n-5}(\cV_l,
\bZ) \oplus H^{n-1}(T_l, \bZ) \oplus
\]
\[
\oplus H^{n-3}(T_l, \bZ) \oplus
H^{n-5}(T_l, \bZ) \oplus H^{n-7}(T_l, \bZ) \oplus H^{n-9}(T_l, \bZ)
\oplus H^{n-11}(T_l, \bZ) \:\: .
\]

\item[$\bold E''$] $\bP^2$-bundle over $R_l'$ which contains $E'' \cap
Q_l''$ as a conic-bundle over $R_l'$. Hence
\[
H^{n-1}(E'', \bZ) \cong H^{n-1}(R_l', \bZ) \oplus H^{n-3}(R_l', \bZ)
\oplus H^{n-5}(R_l', \bZ) \:\: .
\]

\end{description}
\end{lemma}

{\em Proof :} Easy. \hfill \qed

\begin{lemma}
There is a natural exact sequence
\[
0 \lra H^{n-3}(T_l, \bQ) \oplus H^{n-5}(T_l, \bQ) \oplus H^{n-7}(T_l,
\bQ)^{\oplus 2} \oplus H^{n-9}(T_l, \bQ) \oplus H^{n-3}(R_l', \bQ)
\lra
\]
\[
\lra H^{n-1}(Q_l'', \bQ) \lra H^{n+1}(\bP^{n-2}, \bQ) \lra 0
\]
where the map
\[
H^{n-3}(T_l, \bQ) \oplus H^{n-5}(T_l, \bQ) \oplus H^{n-7}(T_l,
\bQ)^{\oplus 2} \oplus H^{n-9}(T_l, \bQ) \lra H^{n-1}(Q_l'', \bQ)
\]
is obtained from the inclusion $T_l'' \subset Q_l''$.
\label{lemHQl}
\end{lemma}

{\em Proof :} From the description of $\tE' \cap Q_l''$ in Lemma
\ref{descrinters} it follows that the map $d_{n-3}$ in sequence
(\ref{Qlcoh}) is injective and we have the exact sequence
\begin{eqnarray}
0 \lra H^{n-5}(\tE' \cap E'' \cap Q_l'', \bZ) \stackrel{d_{n-3}}{\lra}
H^{n-3}(\tE' \cap E'', \bZ) \oplus H^{n-3}(\tE' \cap Q_l'', \bZ)
\oplus H^{n-3}(E'' \cap Q_l'', \bZ) \stackrel{d_{n-1}}{\lra} \nonumber
\end{eqnarray}
\begin{eqnarray}
\stackrel{d_{n-1}}{\lra} H^{n-1}(\tE', \bZ) \oplus H^{n-1}(E'', \bZ)
\oplus H^{n-1}(Q_l'', \bZ) \stackrel{d_{n+1}}{\lra} H^{n+1}({\bP^{n-2}}'',
\bZ) \lra 0 \:\: .
\nonumber
\end{eqnarray}
Tensoring the exact sequence (\ref{Qlcoh}) with $\bQ$ and using Lemma
\ref{descrinters} and the isomorphism
\[
H^{n+1}({\bP^{n-2}}'', \bZ) \cong H^{n+1}(\bP^{n-2}, \bZ) \oplus
\]
\[
\oplus
H^{n-1}(T_l, \bZ) \oplus H^{n-3}(T_l, \bZ) \oplus H^{n-5}(T_l, \bZ)
\oplus H^{n-7}(T_l, \bZ) \oplus H^{n-9}(T_l, \bZ) \oplus
\]
\[
\oplus H^{n-1}(R_l',
\bZ) \oplus H^{n-3}(R_l', \bZ) \:\: ,
\]
we easily deduce Lemma \ref{lemHQl}. \hfill $\qed$

\begin{remark}
In fact we have the exact sequence
\[
0 \lra H^{n-3}(T_l, \bZ [\frac{1}{30}]) \oplus H^{n-5}(T_l, \bZ
[\frac{1}{30}]) \oplus H^{n-7}(T_l, \bZ [\frac{1}{30}])^{\oplus 2}
\oplus H^{n-9}(T_l, \bZ [\frac{1}{30}]) \oplus H^{n-3}(R_l', \bZ
[\frac{1}{30}])
\]
\[
\lra H^{n-1}(Q_l'', \bZ [\frac{1}{30}]) \lra
H^{n+1}(\bP^{n-2}, \bZ [\frac{1}{30}]) \lra 0 \:\: .
\]
\end{remark}

It follows from the previous lemma (since the cohomology of $X$ has no
torsion) that the image of
\[
H^{n-5}(Q_l'', \bZ) \lra H^{n-1}(X, \bZ)
\]
is algebraic. Hence the image of the composition $H^{n-5}(Q_l'', \bZ)
\lra H^{n-1}(X, \bZ) \surj H^{n-1}(X, \bZ)^0$ is algebraic. For $X$
generic $H^{n-1}(X, \bZ)^0$ has no nonzero algebraic part. Hence for
$X$ generic and therefore, for all $X$, the image of $H^{n-5}(Q_l'',
\bZ) \lra H^{n-1}(X, \bZ)^0$ is zero. Hence the map
\[
H^{n-3}(L_l, \bZ) \lra H^{n-1}(X, \bZ)^0
\]
is onto. We have
\[
H^{n-3}(L_l, \bZ) \cong H^{n-3}(S_l', \bZ) \oplus H^{n-5}(S_l', \bZ)
\:
\]
and the restriction $H^{n-5}(S_l', \bZ) \lra H^{n-1}(X, \bZ)^0$ is the
composition of pull-back $H^{n-5}(S_l', \bZ) \lra H^{n-5}(S_l'',
\bZ)$, and push-forward $H^{n-5}(S_l'', \bZ) \lra H^{n-5}(Q_l'', \bZ)
\lra H^{n-1}(X, \bZ)^0$. Hence the map $H^{n-5}(S_l', \bZ) \lra
H^{n-1}(X, \bZ)^0$ is zero and the map
\[
H^{n-3}(S_l', \bZ) \lra H^{n-1}(X, \bZ)^0
\]
is onto.

Now, we have
\[
H^{n-3}(S_l', \bZ) \cong H^{n-3}(S_l, \bZ) \oplus H^{n-5}(\cP^*_l, \bZ)
\oplus H^{n-7}(\cP^*_l, \bZ) \: .
\]
recall that $\cP^*_l$ is the variety parametrizing lines in the fibers
of $\pi^{-1}(T_l) \lra T_l$. Therefore $\cP^*_l$ is a $\bP^2$-bundle
over $T_l$. Using the fact that $T_l$ is a smooth complete
intersection of dimension $n-8$ in $\bP^{n-2}$, one immediately sees
that the image of the summand $H^{n-5}(\cP^*_l, \bZ) \oplus
H^{n-7}(\cP^*_l, \bZ)$ of $H^{n-3}(S_l', \bZ)$ in $H^{n-1}(X, \bZ)^0$
is zero. Therefore the map
\[
H^{n-3}(S_l, \bZ) \lra H^{n-1}(X, \bZ)^0
\]
is onto. This proves the theorem in the case where $n$ is even, since in
that case $H^{n-1}(X, \bZ)^0 = H^{n-1}(X, \bZ)$.

Let $\sigma_1$ be the inverse image in $S_l$ of the hyperplane class
on the Grassmannian $G(2,n+1)$ by the composition $S_l \lra D_l \inj
G(2,n+1)$. If $n$ is odd, one easily computes that the image of
$\sigma_1^{(n-3)/2}$ in $H^{n-1}(X, \bZ)$ is $5 \zeta^{(n-1)/2}$ where
$\zeta$ is the hyperplane class on $X$. On the
other hand, let $x$ be a general point on $l$ and let $L_x$ be the
union of the lines in $X$ through $x$. Then $L_x$ is the intersection
of $X$ with the hyperplane tangent to $X$ at $x$ and a quadric (it is
the second osculating cone to $X$ at $x$). The cohomology class of a
linear section (through $x$) of $L_x$ of codimension $\frac{n-1}{2} -
2$ is $2 \zeta^{(n-1)/2}$ in $X$ and it is in the image of $H^{n-3}(S_l,
\bZ)$. Since $2$ and $5$ are coprime, the image of $H^{n-3}(S_l,
\bZ)$ in $H^{n-1}(X, \bZ)$ contains $\zeta^{(n-1)/2}$ and the map
\[
\psi : H^{n-3}(S_l, \bZ) \lra H^{n-1}(X, \bZ)
\]
in onto for $n$ odd as well. It is now immediate that $\psi'$ is also
onto for $n$ odd. \hfill $\qed$

\vskip20pt

Let $h$ be the first Chern class of the pull-back of
$\cO_{\bP^{n-2}}(1)$ to $S_l'$, let $\sigma_i$ be the pull-back to
$S'_l$ of the $i$-th Chern class of the universal quotient bundle on
the grassmannian $G(2,n+1) \supset D_l$ and let $e_2$ be the first
Chern class of the exceptional divisor of $S'_l \lra S_l$. We make the

\begin{definition}
For a positive integer $k$ the $k-th$ primitive cohomologies of $S_l$ and
$S_l'$ are
\[
H^k(S_l, \bZ)^0 := (\bZ h \oplus \bZ \sigma_1)^{\perp} \subset
H^k(S_l, \bZ)
\]
and
\[
H^k(S_l', \bZ)^0 := (\bZ h \oplus \bZ \sigma_1 \oplus \bZ e_2)^{\perp}
\subset H^k(S_l', \bZ)
\]
where $\perp$ means orthogonal complement with respect to cup-product.
\label{defprim}
\end{definition}

Composing the map $\psi'$ with restriction to $H^{n-3}(S_l', \bZ)^0$ on
the right and with the projection $H^{n-1}(X, \bZ) \surj H^{n-1}(X,
\bZ)^0$ on the left, we get ${\psi'}^0 : H^{n-3}(S_l', \bZ)^0 \lra
H^{n-1}(X, \bZ)^0$. Our goal is to prove the following generalization
of the results of Clemens and Griffiths.

\begin{theorem}
The map ${\psi'}^0$ is onto and its kernel is the $i_l$-invariant part
$H^{n-3}(S_l', \bZ)^{0+}$ of $H^{n-3}(S_l', \bZ)^0$.
\label{maintheorem}
\end{theorem}

The first step for proving the theorem is

\begin{theorem}
Let $a$ and $b$ be two elements of $H^{n-3}(S_l', \bZ)^0$. Then
\[
\psi'(a). \psi'(b) = a. i_l^* b - a . b \: .
\]
\label{thmaibab}
\end{theorem}
{\em Proof :}
We have
\[
\psi'(a). \psi'(b) = (\epsilon_1 \epsilon_2 \rho)_* p^* a . (\epsilon_1
\epsilon_2 \rho)_* p^* b = (\epsilon_2 \rho)_* p^* a. \epsilon_1^*
{\epsilon_1}_* (\epsilon_2 \rho)_* p^* b \: .
\]
Let $\xi_1$ be the first Chern class of the tautological invertible
sheaf for the projective bundle $g_1 : E_1 \lra l$. Let $\gamma_i^1$ be the
Chern classes of the universal quotient bundle on the projective
bundle $g_1 : E_1 \lra l$, i.e.,
\[
\gamma^1_i = \xi_1^i + \xi_1^{i-1}.g_1^* c_1(N_{l/X}) + ... + g_1^*
c_i(N_{l/X}) \: .
\]
Define $\xi_2, \gamma^2_i$ and $\xi_3, \gamma^3_i$ similarly for the
projective bundles $g_2 : E_2 \lra \pi^{-1}(T_l)$ and $g_3 : E_3 \lra
Q_l''$ respectively.
By, e.g., \cite{B2}, 0.1.3, we have
\[
\epsilon_1^*
{\epsilon_1}_* (\epsilon_2 \rho)_* p^* b = (\epsilon_2 \rho)_* p^* b +
{i_1}_* \left( \sum_{r=0}^{n-4} \xi_1^r.g_1^* {g_1}_* \left(
\gamma^1_{n-4-r}.i_1^* \left( (\epsilon_2 \rho)_* p^* b \right)
\right) \right)
\]
where $i_1 : E_1 \inj X_l$ is the inclusion. We also let $i_2 : E_2
\inj X_l'$ and $i_3 : E_3 \inj X_l''$ be the inclusions.

For any $r$, ($0 \leq r \leq n-4$), we have
\[
{g_1}_* \left( \gamma_{n-4-r} . i_1^* (\epsilon_2 \rho)_* p^* b
\right) \in H^{n-3-2r}(l, \bZ) \: .
\]
Therefore ${g_1}_* \left( \gamma_{n-4-r} . i_1^* (\epsilon_2 \rho)_* p^* b
\right) \neq 0$ only if $n-3-2r = 0$ or $n-3-2r = 2$. This is impossible if
$n$ is even so {\em we now suppose that $n$ is odd}. So if we put
\[
B := {i_1}_* \left( \xi_1^{(n-3)/2}.g_1^* {g_1}_* \left(
\gamma^1_{(n-5)/2}.i_1^* \left( (\epsilon_2 \rho)_* p^* b \right)
\right) \hbox{\linebreak}
+ \xi_1^{(n-5)/2}.g_1^* {g_1}_* \left(
\gamma^1_{(n-3)/2}.i_1^* \left( (\epsilon_2 \rho)_* p^* b \right)
\right) \right) \: ,
\]
we have
\[
\epsilon_1^*
{\epsilon_1}_* (\epsilon_2 \rho)_* p^* b = (\epsilon_2 \rho)_* p^* b + B \: .
\]
If $n \geq 7$, replacing $\gamma^1_{(n-5)/2}$ and $\gamma^1_{(n-3)/2}$
in terms of $\xi_1$, we obtain
\[
B = {i_1}_* \left( \xi_1^{(n-3)/2}.g_1^* {g_1}_* \left(
\xi_1^{(n-5)/2} . i_1^* \left( (\epsilon_2 \rho)_* p^* b \right) +
\xi_1^{(n-7)/2} . g_1^* c_1(N_{l/X}) . i_1^* \left( (\epsilon_2
\rho)_* p^* b \right)  \right) \right) +
\]
\[
+ {i_1}_* \left( \xi_1^{(n-5)/2}.g_1^* {g_1}_*
\left( \xi_1^{(n-3)/2} . i_1^* \left( (\epsilon_2 \rho)_* p^* b \right) +
\xi_1^{(n-5)/2} . g_1^* c_1(N_{l/X}) . i_1^* \left( (\epsilon_2
\rho)_* p^* b \right)  \right)
\right) \: .
\]
We have $c_1(N_{l/X}) = (n-4)j_1^* \zeta$ where $\zeta= c_1(\cO_{\bP^n}(1))$
and $j_1 : l \inj X$ is the inclusion. Similarly we define $j_2 :
\pi^{-1}(T_l) \inj X_l$ and $j_3 : Q_l'' \inj X_l'$ to be the
inclusions. Therefore we obtain
\[
B = {i_1}_* \left( \xi_1^{(n-3)/2}.g_1^* {g_1}_* \left(
\xi_1^{(n-5)/2} . i_1^* \left( (\epsilon_2 \rho)_* p^* b \right) +
\xi_1^{(n-7)/2} . (n-4) g_1^* j_1^* \zeta . i_1^* \left( (\epsilon_2
\rho)_* p^* b \right)  \right) \right) +
\]
\[
+ {i_1}_* \left( \xi_1^{(n-5)/2}.g_1^* {g_1}_*
\left( \xi_1^{(n-3)/2} . i_1^* \left( (\epsilon_2 \rho)_* p^* b \right) +
\xi_1^{(n-5)/2} . (n-4) g_1^* j_1^* \zeta . i_1^* \left( (\epsilon_2
\rho)_* p^* b \right)  \right)
\right) \: .
\]
Or, since $j_1 g_1 = \epsilon_1 i_1$,
\[
B = {i_1}_* \left( \xi_1^{(n-3)/2}.g_1^* {g_1}_* \left(
\xi_1^{(n-5)/2} . i_1^* \left( (\epsilon_2 \rho)_* p^* b \right) +
\xi_1^{(n-7)/2} . (n-4) i_1^* \epsilon_1^* \zeta . i_1^* \left( (\epsilon_2
\rho)_* p^* b \right)  \right) \right) +
\]
\[
+ {i_1}_* \left( \xi_1^{(n-5)/2}.g_1^* {g_1}_*
\left( \xi_1^{(n-3)/2} . i_1^* \left( (\epsilon_2 \rho)_* p^* b \right) +
\xi_1^{(n-5)/2} . (n-4) i_1^* \epsilon_1^* \zeta . i_1^* \left( (\epsilon_2
\rho)_* p^* b \right)  \right)
\right) \: .
\]
Let $E_1$ also denote the first Chern class of the invertible sheaf
$\cO_{X_l}(E_1)$. Since $\xi_1 = - i_1^* E_1$, we can write
\[
B =(-1)^{n} {i_1}_* \left(  i_1^* E_1^{(n-3)/2}.g_1^* {g_1}_* i_1^*\left(
 E_1^{(n-5)/2} . \left( (\epsilon_2 \rho)_* p^* b
 \right) - E_1^{(n-7)/2} . (n-4) \epsilon_1^*
 \zeta . \left( (\epsilon_2 \rho)_* p^* b \right)
 \right) \right) +
\]
\[
+ (-1)^n {i_1}_* \left( i_1^* E_1^{(n-5)/2}.g_1^* {g_1}_* i_1^*\left(
 E_1^{(n-3)/2} . \left( (\epsilon_2 \rho)_* p^* b
\right) - E_1^{(n-5)/2} . (n-4) \epsilon_1^*
\zeta . \left( (\epsilon_2 \rho)_* p^* b \right)
\right) \right) \: .
\]
Or, since ${g_1}_*i_1^*={j_1}^*{\epsilon_1}_*$,
\[
B =(-1)^{n} {i_1}_* \left(  i_1^* E_1^{(n-3)/2}.g_1^*
{j_1}^*{\epsilon_1}_* \left(
 E_1^{(n-5)/2} . \left( (\epsilon_2 \rho)_* p^* b
 \right) - E_1^{(n-7)/2} . (n-4) \epsilon_1^*
\zeta . \left( (\epsilon_2 \rho)_* p^* b \right)
 \right) \right) +
\]
\[
+ (-1)^n {i_1}_* \left( i_1^* E_1^{(n-5)/2}.g_1^*
{j_1}^*{\epsilon_1}_* \left(
 E_1^{(n-3)/2} . \left( (\epsilon_2 \rho)_* p^* b
\right) - E_1^{(n-5)/2} . (n-4) \epsilon_1^*
\zeta . \left( (\epsilon_2 \rho)_* p^* b \right)
\right) \right) \: .
\]
Now
\[
{\epsilon_1}_* \left( E_1^{(n-5)/2} . \left( (\epsilon_2 \rho)_* p^* b
 \right) - E_1^{(n-7)/2} . (n-4) \epsilon_1^*
 \zeta . \left( (\epsilon_2 \rho)_* p^* b \right)\right)
\]
is an element of $H^{2n-6}(X, \bZ)$. Hence its image by $j_1^*$ is
zero unless $2n-6 \leq 2$, i.e., $n \leq 4$. We supposed $n \geq 7$.
Similarly,
\[
j_1^*{\epsilon_1}_* \left( E_1^{(n-5)/2} . \left( (\epsilon_2 \rho)_* p^* b
 \right) - E_1^{(n-7)/2} . (n-4) \epsilon_1^*
 \zeta . \left( (\epsilon_2 \rho)_* p^* b \right)\right)
\]
is zero unless $2n-4 \leq 2$ which implies $n \leq 3$. Hence $B$ is
zero for $n \geq 7$. Similarly, $B$ is zero for $n=5$.

\vskip10pt
Therefore
\[
\psi'(a) . \psi'(b) = (\epsilon_2 \rho)_* p^* a. (\epsilon_2 \rho)_* p^*
b \: .
\]

\vskip10pt
Now write
\[
\psi'(a) . \psi'(b) = \rho_* p^* a. \epsilon_2^* {\epsilon_2}_* \rho_*
p^* b
\]
and, as before,
\[
\epsilon_2^* {\epsilon_2}_* \rho_* p^* b = \rho_* p^* b + {i_2}_*
\left( \sum_{r=0}^3 \xi_2^r . g_2^* {g_2}_* \left( \gamma^2_{3-r}
. i_2^* \rho_* p^* b \right) \right) \: .
\]
So
\[
\psi'(a) . \psi'(b) = \rho_* p^* a. \rho_* p^* b + \rho_* p^* a . {i_2}_*
\left( \sum_{r=0}^3 \xi_2^r . g_2^* {g_2}_* \left( \gamma^2_{3-r}
. i_2^* \rho_* p^* b \right) \right)
\]
or
\[
\psi'(a) . \psi'(b) = \rho_* p^* a. \rho_* p^* b + i_2^* \rho_* p^* a .
\left( \sum_{r=0}^3 \xi_2^r . g_2^* {g_2}_* \left( \gamma^2_{3-r}
. i_2^* \rho_* p^* b \right) \right) \: .
\]
We have $a. e_2 = 0$. Hence $p^* a . p^* e_2 = 0$. Let $E_2$ also
denote the cohomology class of $E_2$. Then it is easily seen that
$\rho^* E_2 = p^* e_2$. Therefore $p^* a . \rho^* E_2 = 0$. In order
to use this, we need to modify the above expression a bit.

We first need to write the first three Chern classes of
$N_{\pi^{-1}(T_l)/X_l}$ as inverse images of cohomology classes by
$j_2$. Consider the exact sequence
\[
0 \lra N_{\pi^{-1}(T_l)/X_l} \lra N_{\pi^{-1}(T_l)/ \bP^n_l} \lra
N_{X_l/ \bP^n_l}|_{\pi^{-1}(T_l)} \lra 0 \: .
\]
We have
\[
N_{X_l/ \bP^n_l} \cong \cO_{\bP E}(2) \otimes \pi^*
\cO_{\bP^{n-2}}(1)
\]
where $E = \cO_{\bP^{n-2}}(-1) \oplus \cO_{\bP^{n-2}}^{\oplus 2}$, so that
$\bP E \cong \bP^n_l$. Also
\[
N_{\pi^{-1}(T_l)/ \bP^n_l} \cong \pi^*N_{T_l/ \bP^{n-2}} \cong
\pi^* \left( \cO_{\bP^{n-2}}(3) \oplus \cO_{\bP^{n-2}}(2)^{\oplus 2} \oplus
\cO_{\bP^{n-2}}(1)^{\oplus 3} \right) \: .
\]
It follows that we can write $c_i(N_{\pi^{-1}(T_l)/X_l}) = j_2^* c_i$
where the $c_i$ are cohomology classes on $X_l$. So
\[
\gamma^2_r = \xi_2^r + \xi_2^{r-1} . g_2^* j_2^* c_1 + ... + g_2^*
j_2^* c_r
\]
and, since $\xi_2 = - i_2^* E_2$ and $j_2 g_2 = \epsilon_2 i_2$, we have
\[
\gamma^2_r = i_2^* \alpha^2_r
\]
where
\[
\alpha^2_r = (-1)^r E_2^r + (-1)^{r-1} E_2^{r-1}
. \epsilon_2^* c_1 + ... + \epsilon_2^* c_r \: .
\]
Therefore, using ${g_2}_* i_2^* = j_2^* {\epsilon_2}_*$ and $j_2 g_2 =
\epsilon_2 i_2$,
\[
i_2^* \rho_* p^* a . \left( \sum_{r=0}^3 \xi_2^r . g_2^* {g_2}_*
\left( \gamma^2_{3-r} . i_2^* \rho_* p^* b \right) \right) = i_2^*
\left( \rho_* p^* a . \left( \sum_{r=0}^3 (-1)^r E_2^r . \epsilon_2^*
{\epsilon_2}_* \left( \alpha^2_{3-r} . \rho_* p^* b \right) \right)
\right)
\]
\[
= \rho_* p^* a . E_2 . \left( \sum_{r=0}^3 (-1)^r E_2^r . \epsilon_2^*
{\epsilon_2}_* \left( \alpha^2_{3-r} . \rho_* p^* b \right) \right) =
p^* a . \rho^* E_2 . \rho^* \left( \sum_{r=0}^3 (-1)^r E_2^r . \epsilon_2^*
{\epsilon_2}_* \left( \alpha^2_{3-r} . \rho_* p^* b \right) \right) = 0
\]
and we obtain
\[
\psi'(a) . \psi'(b) = \rho_* p^* a. \rho_* p^* b \: .
\]

Writing $\rho = \epsilon_3 \tq$, we have
\[
\psi'(a) . \psi'(b) = (\epsilon_3 \tq)_* p^* a. (\epsilon_3 \tq)_* p^* b
= \tq_* p^* a . \epsilon_3^* {\epsilon_3}_* \tq_* p^* b
\]
and, as before,
\[
\psi'(a) . \psi'(b) = \tq_* p^* a . \tq_* p^* b + \tq_* p^* a . {i_3}_*
g_3^* {g_3}_* i_3^* \tq_* p^* b = \tq_* p^* a . \tq_* p^* b + i_3^*
\tq_* p^* a . g_3^* {g_3}_* i_3^* \tq_* p^* b \: .
\]
Consider the commutative diagram
\[
\begin{array}{ccccccc}
 & & S_l'' & \stackrel{q'}{\lra} & E_3 & \stackrel{g_3}{\lra} & Q_l'' \\
 & \stackrel{\epsilon_4}{\swarrow} & \downarrow i_3' & & \downarrow i_3 &
 & \downarrow j_3 \\
S_l' & \stackrel{p}{\longleftarrow} & L_l & \stackrel{\tq}{\lra} & X_l''
& \stackrel{\epsilon_3}{\lra} & X_l'
\end{array} 
\]
where the two squares are fiber squares. Using the diagram, we modify
$\psi'(a). \psi'(b)$ as follows
\[
\psi'(a) . \psi'(b) = \tq_* p^* a . \tq_* p^* b + q'_* {i_3'}^* p^* a
. g_3^* {g_3}_* q'_* {i_3'}^* p^* b = 
\]
\[
= \tq_* p^* a . \tq_* p^* b + q'_*
{\epsilon_4}^* a . g_3^* {g_3}_* q'_* {\epsilon_4}^* b = \tq_* p^* a
. \tq_* p^* b + {\epsilon_4}^* a . (g_3 q')^* (g_3 q')_* {\epsilon_4}^* b
\: .
\]
The morphism $g_3q' : S_l'' \lra Q_l''$ is a double cover whose
involution $i_l'$ is the lift of $i_l$. Therefore
\[
(g_3 q')^* (g_3 q')_* {\epsilon_4}^* b = \epsilon_4^* b + {i_l'}^*
\epsilon_4^* b = \epsilon_4^* b + \epsilon_4^* i_l^* b
\]
and
\[
{\epsilon_4}^* a . (g_3 q')^* (g_3 q')_* {\epsilon_4}^* b = \epsilon_4^* a
. \left( \epsilon_4^* b + \epsilon_4^* i_l^* b \right) = a
. {\epsilon_4}_* \left( \epsilon_4^* b + \epsilon_4^* i_l^* b \right)
= a . \left( b + i_l^* b \right) \: .
\]
On the other hand
\[
\tq_* p^* a . \tq_* p^* b = p^* a . p^* b . \tq^* L_l
\]
where we also denote by $L_l$ the cohomology class of $L_l$ in $X_l''$. We
have the following
\begin{lemma}
The cohomology class of $L_l$ in $X_l''$ is equal to
\[
5 \left(\epsilon_1 \epsilon_2 \epsilon_3 \right)^* \zeta
- 5 \left( \epsilon_2 \epsilon_3 \right)^* E_1 - 2 E_3 - k
\epsilon_3^* E_2
\]
for some nonnegative integer $k$.
\end{lemma}
{\em Proof :}
To compute the coefficient of $\left(\epsilon_1 \epsilon_2 \epsilon_3
\right)^* \zeta$, we push $L_l$ forward to $X$ and
compute its degree in $\bP^n$. The image of $L_l$ in $X$ is the union
of all the lines in $X$ which are incident to $l$. Since any such line
maps to a point of $Q_l$ by the projection from $l$, the image of
$L_l$ is the intersection with $X$ of the cone of vertex $l$ over
$Q_l$. Since $Q_l$ has degree $5$, this proves that the coefficient of
$\left(\epsilon_1 \epsilon_2 \epsilon_3 \right)^* \zeta$
is $5$.

The coefficient of $\left( \epsilon_2 \epsilon_3 \right)^* E_1$ is the
negative of the multiplicity of the image of $L_l$ in $X$ along
$l$. Intersecting $X$ with a general linear subspace of dimension $3$
which contains $l$, we see that this linear subspace contains $10$
distinct lines which are distinct from $l$ and are in the image of
$L_l$. Therefore, the multiplicity of the image of $L_l$ along $l$ is
exactly $5 = 5.3 - 10$.

The coefficient of $E_3$ is the negative of the multiplicity of the
image of $L_l$ in $X_l'$ along $Q_l''$. This is $2$ since $L_l$ is
smooth and $\rho$ is an embedding outside $S_l''$ and has degree $2$
on $S_l''$.
\hfill
$\qed$

\vskip30pt
Now we will use the hypothesis $a.h = 0$. It implies $p^* a . p^*h =
0$. One easily sees that
\[
p^* h = \left( \epsilon_2 \rho \right)^* \pi_X^* c_1(\cO_{\bP^{n-2}}(1))
\: .
\]
On the other hand $\epsilon_1^* \zeta - E_1 = \pi^*
c_1(\cO_{\bP^{n-2}}(1))$. Therefore
\[
p^* a . \left(\epsilon_1 \epsilon_2 \rho \right)^*
\zeta = p^* a . \left( \epsilon_2 \rho \right)^*
E_1 \: .
\]
Furthermore, we saw that $p^* a . \rho^* E_2 = 0$, hence,
\[
\tq_* p^* a . \tq_* p^* b = p^* a . p^* b . \tq^* L_l = p^* a . p^* b
. \left( - 2 \tq^* E_3 \right) = - 2 a . b \: \: .
\]
Finally,
\[
\psi'(a) . \psi'(b) = - 2 a . b + a . \left( b + i_l^* b \right) = a
. i_l^* b - a . b \: \: .
\]
\hfill
$\qed$
\vskip.5in
\begin{corollary} \label{kerpsi0}
If ${\psi'}^0$ is onto, the kernel of ${\psi'}^0$ is equal to the set
of $i_l$-invariant elements of $H^{n-3}(S_l', \bZ)$.
\end{corollary}
{\em Proof :} Let $b$ be an element of $H^{n-3}(S_l', \bZ)^0$. Then
${\psi'}^0(b)$ is zero if and only if,
\[
\hbox{for every element }c \hbox{ of } H^{n-1}(X, \bZ)^0, \:\: \psi'(b) . c
= 0 \: \: .
\]
If ${\psi'}^0$ is onto, this is equivalent to
\[
\hbox{for every element }a \hbox{ of }
H^{n-3}(S_l', \bZ)^0, \: \: \psi'(a) . \psi'(b) = 0 \: \: .
\]
By theorem \ref{thmaibab}, this is equivalent to
\[
\hbox{for every element }a \hbox{ of } H^{n-3}(S_l', \bZ)^0, \:\:
a . \left( i_l^* b - b \right) = 0
\]
which is in turn equivalent to
\[
b = i_l^* b \: \: .
\]
\hfill \qed \vskip20pt

We are now ready to prove
\begin{lemma} \label{NSSl}
Suppose $n \geq 6$, then
\[
H^2(S_l, \bQ) = \bQ h \oplus \bQ \sigma_1
\]
\[
H^2(S_l', \bQ) = \bQ h \oplus \bQ \sigma_1 \oplus \bQ e_2
\]
and, if $n=5$, we have the exact sequence
\[
0 \lra H^2(Q_l , \bZ)^0 \lra H^2(S_l, \bZ)^0 \lra H^4(X, \bZ)^0 \lra 0
\]
and
\[
H^2(S_l, \bQ) = H^2(S_l, \bQ)^0 \oplus \bQ h \oplus \bQ \sigma_1
\]
(note that $T_l = \emptyset$ for $n \leq 7$ so that $Q_l = Q_l'$ and $S_l =
S_l'$).
\end{lemma}
{\em Proof :} First suppose $n=5$. Then the direct sum decomposition above
is clear. To prove the exactness of the sequence, note that $H^2(S_l, \bZ)
\lra H^4(X, \bZ)^0$ is onto by Theorem \ref{thmpsionto}. Since $\bZ h
\oplus \bZ \sigma_1$ is algebraic, its image in $H^4(X, \bZ)^0$ is
algebraic. For $X$ generic, the group $H^4(X, \bZ)^0$ has no nonzero
algebraic part. Therefore for $X$ generic and hence for all $X$, the image
of $\bZ h \oplus \bZ \sigma_1$ in $H^4(X, \bZ)^0$ is zero. It follows that
the sequence is exact on the right. The exactness of the rest of the
sequence now follows from Corollary \ref{kerpsi0}.

Now suppose $n \geq 6$. Since $H^2(S_l', \bQ) \cong H^2(S_l, \bQ) \oplus
\bQ e_2$, we only need to compute $H^2(S_l, \bQ)$. Let $H_1$ be a general
hyperplane in $\bP^{n-2}$ and let $H_2$ be its inverse image in
$\bP^n$. The inverse image $S_{l,H}$ of $H_1$ in $S_l$ parametrizes the
lines in the fibers of $X_{l,H} \lra H_1$ where $X_{l,H}$ is the proper
transform of $X_H := X \cap H_2$ in $X_l$. By \cite{goreskymcpherson} pages
23-25, we have $H^2(S_l, \bZ) \cong H^2(S_{l,H}, \bZ)$ for $n \geq 7$ and
$H^2(S_l, \bZ) \inj H^2(S_{l,H}, \bZ)$ for $n=6$. Suppose therefore that
$n=6$. If we choose a general pencil of hyperplanes in $\bP^{n-2}$ of which
$H_1$ is a member, then $H^2(S_l, \bZ)$ maps into the part of $H^2(S_{l,H},
\bZ)$ which is invariant under monodromy. Since $H^4(X_H, \bZ)^0$ has no
nonzero elements invariant under monodromy, we see that $H^2(S_l, \bZ)^0$
lies in $H^2(Q_{l,H}, \bZ)^0$. Since $H^2(Q_{l,H}, \bZ)^0$ has no nonzero
element invariant under monodromy, we have $H^2(S_l, \bZ)^0=0$ and
$H^2(S_l, \bQ) = \bQ h \oplus \bQ \sigma_1$. \hfill \qed \vskip20pt

We will prove Theorem \ref{maintheorem} in conjunction with some results on
the cohomology of $S_l$ and by induction as follows.

\begin{theorem}
\label{thmpsi0onto}
\begin{enumerate}

\item The maps $\psi^0 : H^{n-3}(S_l, \bZ)^0 \lra H^{n-1}(X, \bZ)^0$ and 
${\psi'}^0 : H^{n-3}(S_l', \bZ)^0 \lra H^{n-1}(X, \bZ)^0$ are onto.
The kernel of  ${\psi'}^0$ is the $i_l$-invariant part
$H^{n-3}(S_l', \bZ)^{0+}$ of $H^{n-3}(S_l', \bZ)^0$ and therefore the 
kernel of $\psi^0$ is $H^{n-3}(S_l, \bZ) \cap H^{n-3}(S_l', \bZ)^{0+}$.

\item The cohomology of $S_l$ is torsion in odd degree except in degree
$n-3$.

\item In even degree the rational cohomology of $S_l$ is generated by 
monomials in $h$ and $\sigma_1$ except in degree $n-3$.

\end{enumerate}

\end{theorem}

{\em Proof :} As mentioned above, we proceed by induction on $n$.

We first show that, for any given $n \geq 5$, parts $2$ and $3$ of the
theorem imply part $1$.

Indeed, assume that parts $2$ and $3$ are true for any smooth cubic
hypersurface in $\bP^n$ for a fixed $n$. Let $Sym(h,\sigma_1)$ be the
subvector space of $H^{n-3}(S_l, \bQ)$ generated by monomials in $h$ and
$\sigma_1$ ($Sym(h, \sigma_1)=0$ if $n$ is even). Then, if $n$ is odd, it
follows from numbers $2$ and $3$ that we have the decomposition
\[
H^{n-3}(S_l, \bQ) \cong H^{n-3}(S_l, \bQ)^0 \oplus Sym(h,\sigma_1) \: .
\]

Since $Sym(h,\sigma_1)$ is algebraic, its image in $H^{n-1}(X, \bZ)$ is
also algebraic. For $X$ generic $H^{n-1}(X, \bZ)^0$ has no algebraic
part. Therefore for $X$ generic and hence for all $X$, the image of
$Sym(h,\sigma_1)$ is zero in $H^{n-1}(X, \bZ)^0$. Since the cohomology
of $X$ has no torsion and, by Theorem \ref{thmpsionto}, the map $\psi :
H^{n-3}(S_l, \bZ) \lra H^{n-1}(X, \bZ)$ is onto, it follows that
\[
\psi^0 : H^{n-3}(S_l, \bZ)^0 \lra H^{n-1}(X, \bZ)^0
\]
is onto.

Since $\psi^0$ is the composition of ${\psi'}^0$ with the inclusion 
$H^{n-3}(S_l, \bZ)^0 \inj H^{n-3}(S_l', \bZ)^0$, we deduce that 
${\psi'}^0$ is also onto. The rest of part $1$ is Corollary \ref{kerpsi0}.
\vskip10pt

Now we prove that parts $1$, $2$ and $3$ for $n-1 \geq 5$ imply parts $2$
and $3$ for $n$. Let $H_1, H_2, X_{l,H}, S_{l, H}$ be as in the proof of
Lemma \ref{NSSl}, let $H_1'$ be the proper transform of $H_1$ in
${\bP^{n-2}}'$ and let $X_{l,H}'$ and $S_{l,H}'$ be the proper transforms
of $X_{l,H}$ and $S_{l,H}$ in $X_l'$ and $S_l'$ respectively. By
\cite{goreskymcpherson} pages 23-25, for every $k \leq n-5$, we have
\[
H^k(S_l, \bZ) \cong H^k(S_{l,H}, \bZ)
\]
and
\[
H^{n-4}(S_l, \bZ) \inj H^{n-4}(S_{l,H}, \bZ) \:\: .
\]
In particular, it follows from this and our induction hypothesis that 
$H^{n-3}(S_l, \bQ)$ and $H^{n-4}(S_l, \bQ)$ are the direct sums of their 
primitive parts and their subvector spaces generated by the monomials in 
$h$ and $\sigma_1$. Now it is enough to show that $H^{n-4}(S_l, \bQ)^0=0$.

If we choose a general pencil of hyperplanes in $\bP^{n-2}$ of which $H_1$ 
is a member, then $H^{n-4}(S_l, \bZ)$ maps into the part of 
$H^{n-4}(S_{l,H}, \bZ)$ which is invariant under monodromy. By our 
induction hypothesis, we have the exact sequence
\[
0 \lra H^{n-4}(S_{l,H}, \bZ)^0 \cap H^{n-4}(S_{l,H}', \bZ)^{0+} \lra 
H^{n-4}(S_{l,H}, \bZ)^0 \lra H^{n-2}(X_H, \bZ)^0 \lra 0 \:\: .
\]
Since $H^{n-2}(X_H, \bZ)^0$ has no nonzero elements invariant under 
monodromy, we see that $H^{n-4}(S_l, \bZ)^0$ lies in $H^{n-4}(S_{l,H}, 
\bZ)^0 \cap H^{n-4}(S_{l,H}', \bZ)^{0+}$. Therefore all the elements of 
$H^{n-4}(S_l, \bZ)^0$ are $i_l$-invariant and hence are contained in 
$H^{n-4}(Q_l', \bZ)^0 \subset H^{n-4}(S_l', \bZ)^0$.

Now let
\[
\begin{array}{ccc}
\bP^n & \subset & \bP^{n+1} \\
\downarrow & & \downarrow \\
\bP^{n-2} & \subset & \bP^{n-1}
\end{array}
\]
be a commutative diagram of linear embeddings and projections from $l$. 
Let $Y$ be a general cubic hypersurface in $\bP^{n+1}$ such that $Y \cap 
\bP^n = X$, let $Y_l$ be the blow up of $Y$ along $l$ and let $S_{l,Y}$ be 
the variety parametrizing lines in the fibers of $Y_l \lra \bP^{n-1}$. 
Then, again by \cite{goreskymcpherson} pages 23-25, we have
\[
H^{n-4}(S_l, \bZ) \cong H^{n-4}(S_{l,Y}, \bZ) \:\: .
\]
Let $T_{l,Y}$ be the variety parametrizing the planes in the fibers of 
$Y_l \lra \bP^{n-1}$ and similarly define $Q_{l,Y}$, $Q_{l,Y}'$, $R_{l,Y}'$
and $Q_{l,Y}''$. By Lemma \ref{lemHQl} we have the exact sequence
\[
0 \lra H^{n-2}(T_{l,Y}, \bQ) \oplus H^{n-4}(T_{l,Y}, \bQ) \oplus
H^{n-6}(T_{l,Y}, \bQ)^{\oplus 2} \oplus H^{n-8}(T_{l,Y}, \bQ) \oplus
H^{n-2}(R_{l,Y}',\bQ) \lra
\]
\[
\lra H^n(Q_{l,Y}'', \bQ) \lra H^{n+2}(\bP^{n-1}, \bQ) \lra 0 \:\: .
\]
It is easily seen that the intersection of the subspace
\[
H^{n-2}(T_{l,Y}, \bQ) \oplus H^{n-4}(T_{l,Y}, \bQ) \oplus H^{n-6}(T_{l,Y},
\bQ)^{\oplus 2} \oplus H^{n-8}(T_{l,Y}, \bQ) \oplus H^{n-2}(R_{l,Y}',\bQ)
\]
of $H^n(Q_{l,Y}'', \bQ) \supset H^n(Q_{l,Y}', \bQ)$ with $H^n(S_{l,Y}, \bQ)
\subset H^n(S_{l,Y}'', \bQ)$ is zero. It immediately follows that
$H^{n-4}(S_{l,Y}, \bQ)^0=H^{n-4}(S_l, \bQ)^0=0$.  \vskip10pt

To finish the proof of the theorem all we need to do is to prove the
theorem in the case $n=5$. Suppose therefore that $n=5$. Then part $3$ is
clear. Part $2$ is proved in \cite{voisin} Lemme 3 page 591. Part $1$ is
Lemma \ref{NSSl}.  \hfill $\qed$



\end{document}